\documentclass[10pt]{amsart}
\usepackage[T1]{fontenc}
\usepackage{lmodern}
\usepackage{eucal}

\usepackage{geometry}
\usepackage{fancyhdr}
\usepackage[parfill]{parskip}
\usepackage{lscape}
\usepackage{rotating}

\usepackage{amsmath, amsthm}
\usepackage{amssymb, amsfonts, amscd}
\usepackage{mathrsfs}
\usepackage{xfrac}

\numberwithin{equation}{section}

\usepackage{array}
\usepackage{tabularx}
\usepackage{multirow}
\usepackage{multicol}
\usepackage{booktabs}

\usepackage{xcolor}
\usepackage{float}
\usepackage{caption}
\usepackage{subfig}
\usepackage{mdframed}
\usepackage{mwe}

\usepackage{graphicx}
\usepackage{tikz}
\usepackage{pgf}
\usepackage{tikz-cd}
\usepackage{tikzscale}
\usepackage[all]{xy}

\usetikzlibrary{arrows.meta, positioning}
\usetikzlibrary{
    backgrounds,
    shapes,
    shapes.geometric,
    arrows,
    decorations.markings,
    intersections,
    patterns,
    calc
}

\usepackage{enumitem}
\usepackage{xspace}
\usepackage{comment}
\usepackage{blindtext}
\IfFileExists{algorithm2e.sty}{\usepackage[ruled, vlined]{algorithm2e}}{}

\definecolor{smoked}{RGB}{216, 212, 204}
\definecolor{mauve}{RGB}{200, 55, 171}
\definecolor{apricot}{RGB}{250, 144, 4}
\definecolor{sky}{RGB}{66, 169, 244}
\definecolor{plum}{RGB}{76, 0, 102}
\definecolor{forest}{RGB}{90, 145, 120}
\definecolor{sand}{RGB}{180, 160, 120}
\definecolor{sabbia}{RGB}{207, 185, 151}
\definecolor{aliceblue}{rgb}{0.9, 0.95, 1.0}

\usepackage{hyperref}
\hypersetup{
    colorlinks=true,
    linktoc=all,
    linkcolor=blue,
    citecolor=blue,
    filecolor=blue,
    urlcolor=blue,
    pdftitle={Geometrising discrete representations in PSL(2,R)}
}

\usepackage{cleveref}

\newcommand\C{{\mathbb C}}

\newcommand\HH{{\mathbb H}}
\newcommand\PP{{\mathbb P}}

\newcommand\R{{\mathbb R}}
\newcommand\sph{{\mathbb S}}
\newcommand\Z{{\mathbb Z}}
\newcommand\RP{{\mathbb{RP}}}

\newcommand\pc[1]{{\C\PP^{#1}}}

\newcommand{\pslr}[1]{{\mathrm{PSL}(#1,\R)}}

\DeclareMathOperator{\hol}{hol}

\newcommand*{\defeq}{\mathrel{\vcenter{\baselineskip0.5ex \lineskiplimit0pt
                     \hbox{\scriptsize.}\hbox{\scriptsize.}}}=}

\theoremstyle{plain}
\newtheorem{thm}{Theorem}[section]
\newtheorem{thma}{Theorem}

\newtheorem{lem}[thm]{Lemma}
\newtheorem{prop}[thm]{Proposition}
\newtheorem{cor}[thm]{Corollary}

\newtheorem*{dtconj}{Deroin--Tholozan conjecture}

\theoremstyle{definition}
\newtheorem{defn}[thm]{Definition}

\newtheorem{rmk}[thm]{Remark}

\theoremstyle{remark}

\newtheorem{ques}[thm]{Question}

\tikzstyle{rb} = [rectangle, rounded corners, minimum width=3cm, minimum height=1cm, text width=3cm, text centered, draw=black, fill=blue!30]
\tikzstyle{sb} = [rectangle, minimum width=2cm, minimum height=1cm, text width=3cm, text centered, draw=black, fill=violet!30]
\tikzstyle{co} = [circle, radius=2cm, text width=3cm, text centered, draw=white, fill=orange!50]
\tikzstyle{cr} = [circle, radius=2cm, text width=3cm, text centered, draw=white, fill=red!40]
\tikzstyle{cg} = [circle, radius=2cm, text width=3cm, text centered, draw=white, fill=LimeGreen!40]
\tikzstyle{cb} = [circle, radius=2cm, text width=3cm, text centered, draw=white, fill=blue!30]
\tikzstyle{cv} = [circle, radius=2cm, text width=3cm, text centered, draw=white, fill=violet!30]
\tikzstyle{cy} = [circle, radius=2cm, text width=3cm, text centered, draw=white, fill=yellow!50]
\tikzstyle{cgg} = [circle, radius=2cm, text width=3cm, text centered, draw=white, fill=green!40]
\tikzstyle{cc} = [circle, radius=2cm, text width=3cm, text centered, draw=white, fill=cyan!30]

\title[Geometrising discrete representations in \(\pslr2\)]{Geometrising discrete representations in \(\pslr2\)}

\author{Gianluca Faraco}
\address[Gianluca Faraco]{School of Mathematics, Monash University, Clayton, VIC 3800, Australia}
\email{gianluca.faraco@monash.edu}

\date{\today}

\begin{document}

\subjclass[2020]{}
\dedicatory{}

\begin{abstract}
    We prove the Deroin--Tholozan conjecture for discrete representations into \(\pslr2\). More precisely, we give a complete characterisation of those representations of a closed surface group with discrete image that arise as holonomies of branched hyperbolic structures. The proof combines a factorisation theorem of Rong for maps between Seifert fibred spaces with a local surgery that allows one to reinsert pinched handles while preserving the prescribed holonomy. As a consequence, we also obtain symmetric realisations on finite covers associated with normal torsion-free subgroups of the image.
\end{abstract}

\maketitle
\tableofcontents

\section{Introduction}\label{sec:introduction}

In the present paper, we consider \textit{branched hyperbolic structures} on surfaces. As we shall see in more detail in \S\ref{sec:background}, these structures are geometric structures in the sense of Ehresmann--Thurston \cite{goldman2022geometric, thurston2022geometry}, \textit{i.e.}, geometries on surfaces locally modelled on the hyperbolic plane \(\mathbb H^2\) along with its group of orientation preserving isometries \(\pslr2\). Every such structure naturally determines a homomorphism, say \(\rho \colon \pi_1(S) \longrightarrow \pslr2\), known as the \textit{holonomy representation}, that encapsulates the geometric data of the underlying structure. Although it is an easy matter to determine such a representation, the reverse direction is generally much more challenging. It is therefore natural to ask to what extent a representation encodes such geometric data, specifically determining the conditions under which a representation can be realised as the holonomy of some branched hyperbolic structure on \(S\). The main result of the present paper directly addresses this existence problem for discrete representations.

\begin{thma}\label{thm:discrete-geometrisation}
    Let \(S\) be a closed oriented surface of genus \(g\ge 2\), and let \(\rho\colon\pi_1(S)\longrightarrow\pslr2\) be a representation with discrete image \(G\) and positive volume. Then the orbifold \(\mathcal O=G\backslash\HH^2\) is compact. Moreover, \(\rho\) is geometrisable if and only if either \(\rho\) is faithful or
        \begin{equation}\label{eq:main-volume-condition}
            \textnormal{vol}(\rho)>\textnormal{vol}(\mathcal O),
        \end{equation}
    where \(\textnormal{vol}(\rho)\) is the volume of \(\rho\) and \(\textnormal{vol}(\mathcal O)\) is the hyperbolic area of the orbifold \(\mathcal O\).
\end{thma}

\subsection{Origin of the Deroin--Tholozan conjecture and motivations} Branched hyperbolic structures belong to a broad class of geometric structures on surfaces in which isolated singularities are allowed. Roughly speaking, a branched geometric structure is allowed to have isolated points around which the local charts are finite branched coverings, locally of the form \(z\mapsto z^{k+1}\). Several familiar geometries arise naturally in this way. Translation surfaces, spherical structures and branched hyperbolic structures are induced by singular Riemannian metrics whose singular curvature is concentrated at finitely many cone points. For translation surfaces and branched hyperbolic structures, a branch point of order \(k\) corresponds to a cone angle \(2\pi(k+1)\), and an analogous description holds for spherical structures. Not every branched geometric structure is naturally described in metric terms; branched complex affine structures provide a basic example. All these geometries can be viewed as subgeometries of complex projective geometry locally modelled on \((\pc1,\operatorname{PSL}(2,\C))\).  
Their branched versions therefore fall under the general framework of branched complex projective structures, whose holonomy problem has been extensively studied.

\subsubsection{Branched hyperbolic structures as complex projective structures} For \(\pc1\)-structures, Gallo--Kapovich--Marden proved that every non-elementary representation into \(\operatorname{PSL}(2,\C)\) is the holonomy of a complex projective structure which is either unbranched or has at most one branch point; see \cite{GalloKapovichMarden}. Their proof left a gap in genus \(2\), which was later completed by Le Fils \cite{LF23}. More recently, Le Fils completed the picture for elementary representations by allowing further branching, and more generally characterised holonomy representations with prescribed branching data; see \cite{LeFilsBPS}. Thus, from the point of view of branched complex projective geometry, the qualitative holonomy problem is completely solved: every representation can be realised once sufficient branching is allowed. A more delicate realisation problem arises when the prescribed representation takes values in the transformation group of a proper subgeometry of complex projective geometry: one may then ask whether it can be realised by a branched structure modelled on that subgeometry, rather than merely by a branched complex projective structure. In the case of interest here, this means asking whether a representation into \(\pslr2\) can be realised as the holonomy of a branched hyperbolic structure. 

There is also a complex-analytic perspective on this realisation problem. A complex projective structure determines a Riemann surface structure and, once this is fixed, its holonomy can be interpreted as the monodromy of a Schwarzian equation, or equivalently of an associated second-order linear ODE, see \cite{GalloKapovichMarden}. For punctured surfaces, the corresponding monodromy problem for Schwarzian equations with regular singularities was studied in \cite{FaracoGupta2025}. This framework is directly related to branched hyperbolic structures. Indeed, removing the branch points from such a structure produces an incomplete hyperbolic structure on the resulting punctured surface. Viewed as a complex projective structure, each puncture is a regular singularity with trivial local projective monodromy, and hence an apparent singularity. In both perspectives, by requiring a representation into \(\pslr2\) to arise from a branched hyperbolic structure nevertheless imposes an additional constraint: the corresponding local projective charts must take values in \(\HH^2\). Thus, even when the prescribed monodromy lies in \(\pslr2\), its realisation by a complex projective structure does not automatically yield a hyperbolic one.


\subsubsection{Branched hyperbolic structures as real projective structures} Branched hyperbolic structures are also examples of branched real projective structures, the real counterpart of the aforementioned complex projective structures. Branched real projective structures and their holonomy representations were studied in \cite{FaracoRungi_anosov,faracorungi_geom}. In the Klein model, the hyperbolic plane \(\HH^2\) is identified with the interior of a conic in \(\RP^2\), and its group of orientation-preserving isometries acts through the irreducible embedding \(\pslr2\longrightarrow\operatorname{PSL}(3,\R)\). Thus branched hyperbolic structures may be viewed as branched real projective structures modelled on a proper subgeometry of \(\RP^2\). This leads to the same realisation problem discussed above: when the holonomy of a branched real projective structure takes values in the transformation group of a proper subgeometry, one may ask whether the representation can be realised by a branched structure modelled on that subgeometry. This question was formulated in \cite[Problem~1.2]{faracorungi_geom}. In the present paper, we address it for the hyperbolic subgeometry.

\subsubsection{Goldman's conjecture and its strong reformulation} The main result of the present paper is also motivated by a number of conjectures, going back to Goldman, concerning the holonomy of branched hyperbolic structures. For the sake of clarity, we briefly recall the relevant background and formulate the version that will be used throughout the paper. Let \(S\) be a closed oriented surface of genus \(g\geq2\). The space of conjugacy classes of representations \(\pi_1(S)\longrightarrow\pslr2\) is the character variety \(\mathfrak X(S)\); see \S\ref{sssec:holomap}. Goldman proved that it has \(4g-3\) connected components, distinguished by the Euler number \(2-2g\leq e(\rho)\leq 2g-2\); see \cite{Goldman88}. Each component with non-zero Euler number is a smooth symplectic manifold of real dimension \(6g-6\); see \cite{Goldman1984}. With the convention adopted in this paper, the volume of a representation is related to its Euler number by the identity \(\textnormal{vol}(\rho)=-2\pi e(\rho)\); see Lemma~\ref{lem:volume-euler}.

The two extremal components, corresponding to \(|e(\,\rho\,)|=2g-2\), are the Fuchsian components and are identified with the classical Teichm\"uller spaces parametrising ordinary hyperbolic structures. The first genuinely branched case is obtained by decreasing the absolute value of the Euler number by one, namely \(|e(\,\rho\,)|=2g-3\): a realisation in either of these components must have total branching order one. This case has played a central role in the development of the subject. Goldman conjectured in \cite{Goldman2006} that every representation in these two components is the holonomy of a branched hyperbolic structure. A first substantial step towards the conjecture was provided by Mathews in \cite{Mathews}, who showed that representations with Euler number \(e=\pm(2g-3)\) are holonomies of branched hyperbolic structures under the additional assumption that some simple closed curve has prescribed holonomy. In genus \(2\), his results, combined with those of March\'e--Wolff in \cite{MarcheWolff2016}, imply that almost every representation with Euler number \(e=\pm1\) is geometrisable. The genus \(2\) case has been recently investigated further by the author and Maret in \cite{FaracoMaret2025}, both for pentagon representations, introduced by Le Fils in \cite{LF23}, and for bow-tie representations, in the course of constructing analogues of Fenchel--Nielsen coordinates on an open and dense subset of the space of holonomy representations of branched hyperbolic structures on a closed surface of genus \(2\). More recently, Deroin and Tholozan proved in \cite{DeroinTholozan2026} that every representation with Euler number \(|e(\,\rho\,)|=2g-3\) is the holonomy of a hyperbolic structure with one simple branch point; the corresponding statement for the opposite sign follows by reversing the orientation. Thus, the one-branch-point case of Goldman's geometrisation problem is now completely understood. Finally, while the present work was being completed, in \cite{Lin2026} Lin obtained a further general result in this direction, proving the existence of branched developing maps for a broad class of representations and, in particular, for almost every representation in every component with non-zero Euler number.

Positive volume alone is not sufficient for geometrisability. Indeed, in his paper \cite{Tan1994} Tan constructed representations of positive volume which do not arise as holonomies of branched hyperbolic structures. A notable feature of his examples is that the image is discrete and the representation has non-trivial kernel. This suggested that the obstruction to geometrisation might be confined to the discrete locus. In this direction, Goldman conjectured in \cite[Conjecture~3.9]{Goldman2006} that every representation with positive volume and dense image is geometrisable. Deroin and Tholozan subsequently proposed a stronger statement, which identifies more precisely the expected exceptional discrete representations and, in particular, implies Goldman's conjecture. We shall refer to their formulation as the \emph{Deroin--Tholozan conjecture}.

\begin{dtconj}\label{conj:deroin-tholozan}
    Let \(\rho\colon\pi_1(S)\longrightarrow\pslr2\) be a representation with positive volume. Then \(\rho\) is the holonomy of a branched hyperbolic structure unless it is non-faithful, has discrete image, and satisfies
    \begin{equation}\label{eq:DT-exception}
        \textnormal{vol}(\rho)
        =
        \textnormal{vol}\bigl(\rho(\pi_1(S))\backslash\HH^2\bigr).
    \end{equation}
\end{dtconj}

The conjecture suggests a striking dichotomy. Representations of positive volume are expected to be geometrisable in great generality, while the exceptional behaviour is concentrated among non-faithful representations with discrete image. In particular, the discrete case is not merely a special instance of the conjecture: it contains the obstruction explicitly singled out by its statement and therefore deserves a separate analysis.

\subsubsection{Discrete representations and mapping class group orbits}

The discrete case already exhibits two rather different behaviours. We shall first focus on discrete representations with cocompact image. Indeed, as we shall see, a discrete representation whose image is not cocompact has Euler number zero and therefore cannot arise as the holonomy of a branched hyperbolic structure, see Lemma \ref{lem:noncompact-euler-zero}. When the discrete image is torsion-free and cocompact, then every non-identity element is hyperbolic and the quotient of \(\HH^2\) is a smooth hyperbolic surface. In this case we shall say that a representation is \textit{purely hyperbolic}. These kind of representations have been already characterised in \cite{Faraco2021}, which settled the torsion-free discrete case under the relevant non-zero degree hypothesis. The purpose of the present paper is to extend this picture to arbitrary discrete representations, including cocompact groups with torsion, and thereby prove the Deroin--Tholozan conjecture for representations with discrete image. There is a geometric reason why the discrete case is particularly well suited to a separate treatment. If \(G=\rho(\pi_1(S))\) is discrete, then \(\mathcal O=G\backslash\HH^2\) is a hyperbolic orbifold, and a developing map can be compared with the map \(S\longrightarrow\mathcal O\). This introduces global topological information which has no direct analogue when the image is dense or otherwise non-discrete. 

At the same time, discrete representations should not be viewed as isolated from the general problem. In every non-maximal component with non-zero Euler number, they form a closed nowhere dense subset of the representation space; see \cite{FunarWolff2007}. Nevertheless, they provide natural reference points from which one can attempt to reach non-discrete representations by deformation. Indeed, a version of the Ehresmann--Thurston principle holds for branched hyperbolic structures: the holonomy of a branched hyperbolic structure remains geometrisable under sufficiently small deformations. Consequently, once the discrete geometrisable representations are understood, it is natural to ask whether a given non-discrete representation can be moved, by changing its marking, arbitrarily close to a discrete representation.

\begin{ques}\label{ques:orbit-discrete}
    Let \(\rho\colon\pi_1(S)\longrightarrow\pslr2\) be a non-discrete representation. Under which conditions does the \(\operatorname{Mod}(S)\)-orbit of \([\rho]\) accumulate on a discrete representation?
\end{ques}

Whenever the \(\operatorname{Mod}(S)\)-orbit of \([\rho]\) accumulates on a discrete geometrisable representation, this has an immediate consequence for the geometrisation problem. Indeed, if a mapping-class translate of the representation \([\rho]\) lies sufficiently close to such a representation, then the Ehresmann--Thurston principle implies that this translate is itself geometrisable. Since precomposition by a mapping class merely changes the marking, the original representation is geometrisable as well.

Question~\ref{ques:orbit-discrete} is naturally related to Goldman's conjecture on the ergodicity of the mapping class group action on the non-maximal components of the \(\pslr2\)-character variety. In genus two, the ergodicity on the components of Euler number \(\pm1\) by March\'e--Wolff \cite{MarcheWolff2016,MarcheWolff2019} whereas the component of representations with null Euler number exhibits a different behaviour. In his recent paper \cite{Lin2026}, Lin showed the ergocidity of the mapping class group in the zero component for \(g\ge3\). Ergodicity would provide a much stronger dynamical statement and, in suitable components, would force generic orbits to visit every positive-measure neighbourhood of the discrete locus. For the purpose of geometrisation, however, such a global measure-theoretic statement is stronger than necessary. Question~\ref{ques:orbit-discrete} asks only for the orbit-closure information needed to approach a geometrisable discrete representation. It therefore isolates a weaker dynamical problem which may already be sufficient for extending the discrete results of this paper to non-discrete representations.

\subsubsection{Hyperbolic analogue of Haupt's Theorem} There is a useful analogy with the period realisation problem for holomorphic differentials. Haupt's theorem characterises the complex-valued characters of \(H_1(S;\Z)\) which arise as periods of holomorphic \(1\)-forms, and Theorem~\ref{thm:discrete-geometrisation} may be viewed as a hyperbolic analogue of its discrete case. To make this comparison explicit, let \(\chi\colon H_1(S;\Z)\longrightarrow\C\) have lattice image \(\Lambda\) and positive volume, where
\begin{equation}\label{eq:flat-volume}
    \textnormal{vol}(\,\chi\,)=\sum_{i=1}^g\operatorname{Im}\bigl(\overline{\chi(a_i)}\,\chi(b_i)\bigr)
\end{equation}
for a symplectic basis \(\{a_1,b_1,\ldots,a_g,b_g\}\). Haupt's theorem states that \(\chi\) is the period character of a holomorphic differential, for some complex structure on \(S\), if and only if
\begin{equation}\label{eq:haupt-condition}
    \textnormal{vol}(\,\chi\,)>\textnormal{vol}(\C/\Lambda);
\end{equation}
see \cite{Haupt1920,Kapovich2020}. The ratio of these two volumes is the degree of the induced holomorphic map to the flat torus \(\C/\Lambda\). Degree one is excluded when \(g\geq2\), since such a map would necessarily be a biholomorphism onto a torus. Theorem~\ref{thm:discrete-geometrisation} yields the analogous degree-one obstruction in the hyperbolic setting, with the quotient hyperbolic orbifold replacing the flat torus. 
This parallel becomes even closer when branching data are prescribed. Refinements of Haupt's theorem due independently to Bainbridge--Johnson--Judge--Park \cite{BJJP} and Le Fils \cite{LeFils} determine the image of the period map on each stratum of holomorphic differentials. In the present paper, we obtain analogous stratum-wise results for branched hyperbolic structures with discrete holonomy, thereby making progress on \cite[Question 1.3]{faracorungi_geom}. As we shall see, the discrete hyperbolic problem exhibits a similar degree obstruction, with the hyperbolic orbifold \(G\backslash\HH^2\) playing the role of the flat torus.

Although these results are motivated by the refined versions of Haupt's theorem, the approach is substantially different. In the presence of torsion, prescribing the branching data leads naturally to an equivariance problem on a finite torsion-free cover. This requires us to study branched hyperbolic structures carrying prescribed finite groups of symmetries. The analysis of such symmetric structures therefore becomes an essential ingredient in the refinement of the geometrisation problem by strata.

\subsection{Symmetric branched hyperbolic structures}\label{ssec:intro symmetries} A first consequence of our main theorem concerns branched hyperbolic structures with prescribed symmetries. Suppose that a branched hyperbolic structure on \(S\) admits a non-trivial finite group \(Q\) of orientation-preserving isometries. The structure then descends to the quotient \(S/Q\), which is a smooth surface when the action is free and, more generally, a hyperbolic cone-metric. Correspondingly, the holonomy representation factors through the fundamental group of the quotient, or through its orbifold fundamental group when fixed points are present. This naturally leads to the problem of constructing branched hyperbolic structures which are invariant under a prescribed finite group action. Questions of this type have appeared previously for projective and branched geometric structures. In \cite{FaracoRuffoni2019}, projective structures with large automorphism groups were studied in relation to highly symmetric Riemann surfaces, while \cite{Faraco2024Automorphisms} investigates automorphism groups of several classes of branched geometric structures and raises related realisation problems with prescribed symmetries.

In the discrete setting considered here, such symmetries arise naturally by passing to finite-index subgroups. Let \(G<\pslr2\) be a discrete cocompact group and let \(H\trianglelefteq G\) be a torsion-free subgroup of finite index. If \(\rho\colon\pi_1(S)\longrightarrow G\) is a representation, then \(\rho^{-1}(\,H\,)\) determines a regular finite covering \(p_H\colon S_H\to S\) whose deck group is naturally isomorphic to \(G/H\). The main theorem implies that, whenever the restricted representation is geometrisable, it admits a realisation carrying this prescribed deck symmetry.

\begin{thma}\label{thma:invariant-cover}
    Let \(\rho\colon\pi_1(S)\to\pslr2\) have discrete cocompact image \(G\). Let \(H\trianglelefteq G\) be torsion-free and of finite index, and let \(p_H\colon S_H\longrightarrow S\) be the covering corresponding to \(\rho^{-1}(\,H\,)\), with the lifted orientation. Suppose \(\rho_H=\rho|_{\pi_1(\,S_H\,)}\) is geometrisable, then it admits a realisation invariant under the prescribed deck action of \(G/H\) on \(S_H\). This realisation is the pull-back of a branched hyperbolic structure on \(S\) with holonomy \(\rho\). In particular, \(G/H\) acts by orientation-preserving isometries of this branched hyperbolic structure.
\end{thma}

We shall also describe when the prescribed \(G/H\)-action coincides with the full group of orientation-preserving isometries of the resulting branched hyperbolic structure. This analysis will later play an essential role in refining the geometrisation problem by prescribing the branching stratum in the presence of torsion.

\smallskip

\subsection{Strata of branched hyperbolic structures}\label{sssec:strata}

The deformation space of branched hyperbolic structures is naturally stratified into strata \(\mathcal{BH}(\,\mu\,)\) enumerated by unordered partitions of positive integers \(0<b<2g-2\). Notice that the case \(b=0\) corresponds to Teichm\"uller space, whereas the Gauss--Bonnet Theorem implies that \(b\leq 2g-3\). For a partition \(\mu=(k_1,\ldots,k_n)\) of some positive integer \(0<b<2g-2\), the stratum \(\mathcal{BH}(\,\mu\,)\) parametrises all branched hyperbolic structures with exactly \(n\) conical singularities of orders \(k_1,\ldots,k_n\). In \S\ref{sec:strata} we consider the refined geometrisation problem in which the branching data are prescribed.

For representations with discrete image, two necessary numerical conditions arise naturally. The Euler number determines the total branching order, while the degree of the induced map to the quotient hyperbolic orbifold bounds the order of each individual branch point. Thus, if \(\rho\) arises as the holonomy of some branched hyperbolic structure in \(\mathcal{BH}(\,\mu\,)\), then
\begin{equation}
    \sum_{i=1}^n k_i=b=e(\,\rho\,)-\chi(S),
\end{equation}
and
\begin{equation}
    k_i+1\leq d(\,\rho\,)=\frac{\textnormal{vol}(\,\rho\,)}{\textnormal{vol}\big(\,G\backslash\HH^2\,\big)}
    \qquad
    \textnormal{for every }i=1,\ldots,n, \textnormal{ where } \rho(\pi_1(S))=G.
\end{equation}

When the image of \(\rho\) is torsion-free, these conditions are also sufficient. We therefore obtain a complete characterisation of the image of the holonomy map on each stratum. 
\begin{thma}\label{thm:intro-strata-torsion-free}
    Let \(\rho\colon\pi_1(S)\to\pslr2\) be a purely hyperbolic representation with positive volume, and let \(\mu=(k_1,\ldots,k_n)\) be a non-empty partition of an integer \(b\), with \(0<b<2g-2\). Then \(\rho\) arises as the holonomy of some branched hyperbolic structure in \(\mathcal{BH}(\,\mu\,)\) if and only if \(e(\,\rho\,)=\chi(S)+b\) and \(k_i+1\leq d(\,\rho\,)\) for every \(i=1,\ldots,n\).
\end{thma}

The proof relies on the realisation theory of primitive branched coverings together with the strong-equivalence theorem of Gabai--Kazez, which allows us to recover the prescribed marking after constructing a covering with the required ramification profile.

When the image contains torsion, the same numerical conditions remain necessary, but an additional equivariance condition appears. In fact, passing to a normal torsion-free finite-index subgroup \(H\trianglelefteq G\) produces a finite cover \(S_H\to S\), on which the restricted representation satisfies the torsion-free criterion. The remaining issue is whether the corresponding realisation can be chosen invariant under the prescribed action of \(G/H\). The results on symmetric branched hyperbolic structures established in \S\ref{sec:sym structures} provide precisely the framework needed to formulate this condition. Thus, the additional condition in the presence of torsion is not numerical, but is encoded by the compatibility between the lifted structure and the prescribed deck action.

\begin{thma}\label{thm:intro-strata-torsion}
    Let \(\rho\colon\pi_1(S)\to\pslr2\) have discrete cocompact image \(G\) with torsion and positive volume, and let \(\mu=(k_1,\ldots,k_n)\) be a non-empty partition of an integer \(b\), with \(0<b<2g-2\). Then \(\rho\) arises as the holonomy of some branched hyperbolic structure in \(\mathcal{BH}(\,\mu\,)\) if and only if
    \begin{equation}
        e(\rho)=\chi(S)+b,
        \qquad
        k_i+1\leq d(\,\rho\,)
        \quad
        \textnormal{for every }i=1,\ldots,n,
    \end{equation}
    and, after passing to a normal torsion-free finite-index subgroup \(H\trianglelefteq G\), the restricted representation admits a realisation in the lifted stratum \(\mathcal{BH}(\,\mu^N\,)\) which is invariant under the prescribed deck action of \(G/H\), where \(N=[G:H]\) and \(\mu^N=(k_1^N,\ldots,k_n^N)\).
\end{thma}

The last condition isolates the additional difficulty created by torsion. Theorem~\ref{cor:invariant-cover} guarantees the existence of a realisation of the restricted representation which is invariant under the prescribed deck action, while Theorem~\ref{thm:intro-strata-torsion-free} guarantees the existence of a realisation in the prescribed lifted stratum. What is required here is that these two properties be achieved simultaneously. Equivalently, the problem reduces to the existence of a symmetric realisation within a prescribed stratum.

\subsection{Idea of the proof}\label{subsec:idea-proof}

The proof of Theorem~\ref{thm:discrete-geometrisation} follows two different routes according to whether the discrete image is torsion-free. Suppose first that \(G=\rho(\pi_1(S))\) is torsion-free, so that \(\Sigma=G\backslash\HH^2\) is a closed hyperbolic surface. Since the surface \(\Sigma\) is aspherical, the representation \(\rho\) is induced by a map \(f\colon S\longrightarrow\Sigma\). We then use Edmonds' factorisation theorem \cite{Edmonds1979}, which deforms a map of non-zero degree between closed surfaces into a composition of surface pinches followed by a branched covering. This reduces the geometrisation problem to the existence criterion established in \cite{Faraco2021} for purely hyperbolic representations. 

This argument, however, does not extend directly when \(G\) contains torsion. In this case \(\mathcal O=G\backslash\HH^2\) is a hyperbolic orbifold rather than a smooth surface, and a map to the underlying topological surface \(|\,\mathcal O\,|\) does not retain the orbifold data carried by the elliptic elements of \(G\). We therefore replace the two-dimensional problem by a three-dimensional one. By starting from a smooth \(\rho\)-equivariant map \(u\colon\widetilde S\to\HH^2\), we realise a principal circle bundle, say \(M\to S\) and the Seifert fibred manifold \(N=G\backslash T^1\HH^2\longrightarrow\mathcal O\). The \(\rho\)-equivariant map \(u\) induces a natural fibre-preserving map, say \(F\colon M\to N\) of positive degree. In this framework, Rong's factorisation theorem \cite{Rong1993} applies and it plays exactly the same role played by Edmonds' theorem for maps of surfaces. It deforms the mapping \(F\) into a composition 
\[ M\xrightarrow{\ \Pi\ }P\xrightarrow{\ B\ }N,
\]  
where \(\Pi\) is a composition of vertical pinches and \(B\) is a fibre-preserving branched covering, see \S\ref{subsec:pinches} for details. The remainder of the proof consists in passing from this factorisation of Seifert manifolds to their base spaces. We show that the intermediate manifold \(P\) is again a principal circle bundle over a smooth surface \(T\), that the vertical pinches induce ordinary pinches \(p\colon S\to T\) of surfaces, and that \(B\) induces a branched map \(T\to|\,\mathcal O\,|\). The pull-back of the orbifold metric has integral cone angles and defines a branched hyperbolic structure on \(T\) whose marked holonomy is the induced representation \(\beta\), with
\[
    \rho=\beta\circ p_*.
\]
Finally, a local handle-insertion surgery reverses the surface pinches while preserving the prescribed holonomy, producing the required branched hyperbolic structure on \(S\).

\subsection{Structure of the paper} The present paper is organised as follows. \S\ref{sec:background} recalls the basic material on branched hyperbolic structures, developing maps and holonomy, Euler number and volume. In \S\ref{sec:seifert-fibrations}, we collect the facts about Seifert fibrations and fibre-preserving maps that will be used throughout the proof of the main theorem, with particular emphasis on exceptional fibres and their local behaviour under branched coverings. In \S\ref{sec:discrete-geometrisation} we prove the main geometrisation theorem for discrete representations. The proof associates to the representation a map between Seifert fibred \(3\)-manifolds, applies Rong's factorisation theorem, and then recovers a branched hyperbolic structure while keeping track of the prescribed holonomy. \S\ref{sec:sym structures} is devoted to symmetric realisations: we apply the main theorem to regular finite coverings associated with normal torsion-free subgroups and study the resulting groups of orientation-preserving isometries. Finally, \S\ref{sec:strata} considers the refined realisation problem with prescribed branching data. We first derive the necessary numerical conditions for a discrete representation to belong to a given stratum, prove that they are sufficient when the image is torsion-free, and then treat the case with torsion by reducing the remaining problem to an equivariant realisation problem on a finite torsion-free cover.

\subsection{Acknowledgements} The author wishes to thank Bertrand Deroin for sharing a preliminary version of his recent paper with Nicholas Tholozan and for asking him about the geometrisation of discrete representations in \(\pslr2\), a question that led to the present work.

\subsection{AI Declaration} All mathematical results in this paper have been developed by the author. The author used ChatGPT-6 Sol to proofread the manuscript and to help rephrase selected passages for clarity and flow. The author takes full responsibility for the content of the paper.

\medskip

\section{Hyperbolic structures}\label{sec:background} Throughout this paper, \(S\) denotes a closed, connected and oriented surface of genus \(g\ge 2\), and let \(\widetilde S\) denote its universal cover. Surfaces are considered with their topological marking and orientation. We choose auxiliary smooth structures for the bundle and differential-form constructions; no complex structure is fixed in advance. We write \(\Gamma=\pi_1(S)\), after choosing a base point, and regard \(\Gamma\) as the group of deck transformations of \(\widetilde S\). The hyperbolic plane \(\HH^2\) is equipped with its standard orientation and metric \(h_{\HH}\) of constant curvature \(-1\). Its group of orientation-preserving isometries is identified with \(\pslr2\) whose action is given by Möbius transformations: an element represented by 
\[
    \pslr2 \times \HH^2 \to \HH^2, \quad \left(\begin{array}{cc} 
a & b\\ c & d \\ \end{array} \right), z \longmapsto \dfrac{az+b}{cz+d}
\]

We shall begin by recalling the geometric structures considered in this paper and the invariants of their holonomy representations. The conventions below are compatible with \cite[\S 2]{Faraco2021}; we fix the sign of the Euler number explicitly.

\subsection{Hyperbolic cone-structures}\label{subsec:bhs}

We begin by recalling the geometric structures considered throughout the paper. For our purposes, we only need branched hyperbolic structures in dimension two, although the definition admits natural generalisations to higher dimensions and to other geometries; see \cite{CHK} for further details.

\begin{defn}[Hyperbolic cone-structure]\label{def: hyp cone-structures}
    A \emph{hyperbolic cone-structure} \(\sigma\) on a surface \(S\) is the datum of a triangulation of \(S\) and a metric such that:
    \begin{itemize}
        \item the link of each vertex is piecewise linearly homeomorphic to a circle;
        \item the restriction of the metric to each simplex is isometric to a geodesic simplex in \(\HH^2\).
    \end{itemize}
\end{defn}

Thus, a two-dimensional hyperbolic cone-structure is obtained by gluing geodesic triangles in \(\HH^2\). The definition also applies to open surfaces and surfaces with possibly geodesic boundary, although in this paper we shall only consider closed surfaces.

Every interior point \(p\in S\) has a neighbourhood locally isometric to \(\HH^2\), except possibly at some vertices of the triangulation, around which the angles sum to a value \(\theta\neq2\pi\). Such points are called \emph{cone points}. A neighbourhood of a cone point is isometric to a wedge of angle \(\theta\) in \(\HH^2\), with its two boundary rays identified. The value \(\theta\) is called the \emph{cone angle}. By setting \(\theta=2\pi(k+1)\), we call \(k\) the \emph{order} of the cone point and write \(\operatorname{ord}(p)=k\). Cone points are the singular points of the cone-structure, whereas all other points are called \emph{regular}. In general, cone angles may be arbitrary positive real numbers and, in particular, may be larger than \(2\pi\). We now introduce two special classes of hyperbolic cone-structures. 

\subsubsection{Branched hyperbolic structures}\label{sssec:bhs} According to Definition \ref{def: hyp cone-structures}, a \emph{branched hyperbolic structure} is a hyperbolic cone-structure for which the order of every cone point is a positive integer. It is an easy matter to check that this is equivalent to the following, generally more common, definition.

\begin{defn}\label{def:branched hy structure}
    A \emph{branched hyperbolic structure} on \(S\), denoted by \(\sigma\), is the datum of a maximal atlas, say \(\mathcal U\), of possibly branched hyperbolic charts. More precisely, a branched hyperbolic structure consists of the following data:
    \begin{itemize}
        \item[1.] a discrete collection of points, say \(\mathcal B\subset S\);
        \item[2.] a maximal atlas of \(\HH^2\)-valued hyperbolic charts on \(S\setminus\mathcal B\); and
        \item[3.] a maximal collection of branched hyperbolic charts centred at the points of \(\mathcal B\),
    \end{itemize}
    such that all transition maps are restrictions of elements of \(\pslr2\). More precisely, for every pair of possibly branched charts \((U_1,\varphi_1)\) and \((U_2,\varphi_2)\) in \(\mathcal U\), and for every connected component, say \(C\subset U_1\cap U_2\), there exists \(g\in\pslr2\) such that \(\varphi_2|_C=g\circ\varphi_1|_C\). At each point \(p\in\mathcal B\), the corresponding branched chart is locally equivalent, in oriented coordinates, to \(z\mapsto z^{k+1}\) for some integer \(k\geq1\).
\end{defn}

We observe that a branched hyperbolic structure naturally determines an underlying complex structure on \(S\). Away from the branch points, the structure is described by orientation-preserving local isometries into \(\HH^2\), and the transition maps belong to \(\pslr2\); hence they are holomorphic for the standard complex structure on \(\HH^2\). Near a branch point of order \(k\), a branched chart is locally of the form \(z\mapsto z^{k+1}\) in an oriented local coordinate \(z\) on the source and a complex coordinate on a disc in \(\HH^2\). The complex structure defined away from the branch points therefore extends across them through the coordinate \(z\). We shall henceforth freely regard a branched hyperbolic structure as a Riemann surface equipped with such local branched hyperbolic charts.

Since \(S\) is closed, hence compact, the set of cone points, also called \textit{branched points} in this setting, is finite. Let \(p_1,\ldots,p_s\) be the branch points, with respective branching orders \(k_1,\ldots,k_s\). We denote by
\begin{equation}
    D_\sigma=\sum_{i=1}^s k_i p_i \qquad \textnormal{ and } \qquad \kappa(\sigma)=\deg D_\sigma=\sum_{i=1}^s k_i
\end{equation}
 the branching divisor of \(\sigma\) and its total branching order respectively. The empty branch set is allowed, in which case \(\sigma\) is an ordinary hyperbolic structure. Cone points may be regarded as points at which curvature is concentrated, and the topology of \(S\) imposes restrictions on the possible total branching. As a consequence of Gauss--Bonnet formula for a cone metric, the identity
\begin{equation}\label{eq:gauss-bonnet}
    -\textnormal{Area}(S,\sigma)+\sum_{i=1}^s\bigl(2\pi-2\pi(k_i+1)\bigr)=2\pi\chi(S)
\end{equation}

holds. Consequently,
\begin{equation}\label{eq:area-branching}
    \textnormal{Area}(S,\sigma)=-2\pi\bigl(\chi(S)+\kappa(\sigma)\bigr)>0.
\end{equation}
In particular, \(\kappa(\sigma)<2g-2\). Thus the total branching order is necessarily bounded by the topology of \(S\).

\smallskip

\subsubsection{Hyperbolic orbifolds}\label{sssec:hyporbi} In the following we shall also make use of hyperbolic orbifolds. These structures are defined as follows.

\begin{defn}[Hyperbolic orbifold]
    An \emph{oriented hyperbolic orbifold} \(\mathcal O\) is an oriented two-dimensional orbifold equipped with a hyperbolic metric. Its underlying topological surface is denoted by \(|\mathcal O|\). Every regular point has a neighbourhood locally isometric to \(\HH^2\), whereas a cone point of order \(m\geq2\) has a neighbourhood modelled on the quotient of a disc in \(\HH^2\) by a rotation of order \(m\). In particular, its cone angle is \(2\pi/m\).
\end{defn}

If \(G<\pslr2\) is discrete, then \(G\backslash\HH^2\) is a hyperbolic orbifold, whose cone points correspond to elliptic elements of \(G\). In particular, if \(G\) is cocompact and torsion-free, then \(G\backslash\HH^2\) is a genuine hyperbolic surface. For a compact oriented hyperbolic orbifold \(\mathcal O\), we shall denote its hyperbolic area by \(\textnormal{vol}(\mathcal O)\). If the underlying surface has genus \(h\) and the cone points have orders \(m_1,\ldots,m_t\), then
\begin{equation}\label{eq:orbifold-area}
    \begin{aligned}
        \chi_{\mathrm{orb}}(\mathcal O) &   = 2-2h-\sum_{i=1}^t\left(1-\frac{1}{m_i}\right),\\ 
                                        &                                                   \\
        \textnormal{vol}(\mathcal O)    &   =  -2\pi\chi_{\mathrm{orb}}(\mathcal O)>0.
    \end{aligned}
\end{equation}

\smallskip

\subsubsection{Quotients by discrete groups}\label{subsubsec:discrete-hyperbolic-quotients} In this section, we describe the local geometry of quotients of the hyperbolic plane by discrete groups and explain why their cone angles are precisely those allowed in a hyperbolic orbifold. We recall the following classical result and include its proof for the reader's convenience.

\begin{prop}\label{prop:discrete-hyperbolic-quotients}
    Let \(G<\mathrm{PSL}(2,\mathbb R)\) be a discrete subgroup. Then the quotient \(\mathcal O=G\backslash\mathbb H^2\), equipped with the quotient metric, is either smooth or an oriented hyperbolic orbifold. More precisely, for each \(p\in\mathbb H^2\), the stabiliser \(G_p\) is cyclic of some finite order \(m\geq 1\), and the cone angle at the image of \(p\) is \(2\pi/m\). In particular, \(\mathcal O\) is a smooth hyperbolic surface if and only if \(G\) is torsion-free. If \(G\) is cocompact, then \(\mathcal O\) is compact.
\end{prop}

\begin{proof}
Fix \(p\in\mathbb H^2\). The stabiliser of \(p\) in \(\mathrm{PSL}(2,\mathbb R)\) is the compact group of rotations about \(p\). Its intersection with the discrete subgroup \(G\) is therefore finite. Every finite subgroup of the rotation group is cyclic, so \(G_p\) is generated by the rotation through \(2\pi/m\), where \(m=|G_p|\). We next justify that this stabiliser determines the local quotient. For every \(R>0\), the set
\begin{equation}
    \bigl\{g\in\mathrm{PSL}(2,\mathbb R): d_{\mathbb H^2}(p,gp)\leq R\bigr\}
\end{equation}
is compact: an orientation-preserving isometry is determined by the image of \(p\) and the image of a unit tangent vector at \(p\), and the unit tangent bundle over a closed hyperbolic ball is compact. Since a discrete subgroup is closed, this set meets \(G\) in finitely many elements. Consequently, there exists \(\varepsilon>0\) such that
\begin{equation}
    d_{\mathbb H^2}(p,gp)>2\varepsilon \qquad\text{for every }g\in G\setminus G_p.
\end{equation}
Thus a neighbourhood of the image of \(p\) in \(\mathcal O\) is isometric to \(B_\varepsilon(p)/G_p\). In geodesic polar coordinates centred at \(p\), the hyperbolic metric is
\begin{equation}
    ds^2=dr^2+\sinh^2(r)\,d\theta^2, \qquad \theta\in\mathbb R/2\pi\mathbb Z.
\end{equation}
Passing to the quotient by \(G_p\) identifies \(\theta\) with \(\theta+2\pi/m\). Hence \(B_\varepsilon(p)/G_p\) is a hyperbolic cone of angle \(2\pi/m\), with \(m=1\) corresponding to a regular point. These local models give the asserted orbifold structure. Finally, since every non-trivial finite-order element of \(\mathrm{PSL}(2,\mathbb R)\) fixes a point of \(\mathbb H^2\), so the quotient is smooth exactly when \(G\) is torsion-free. Compactness readily follows from co-compactness.
\end{proof}

In particular, no cone angle of the quotient metric can exceed \(2\pi\). A hyperbolic cone-structure with an angle greater than \(2\pi\) is therefore not isometric to a quotient \(G\backslash\mathbb H^2\) endowed with its quotient metric.

\subsection{Developing--holonomy pair}\label{subsec:developing-holonomy}

Let \(\sigma\) be a branched hyperbolic structure on \(S\). Analytic continuation of its local charts gives a developing map and a holonomy representation,
\begin{equation}\label{eq:developing-pair}
    \textnormal{dev}\colon\widetilde S\longrightarrow\HH^2, \qquad \rho\colon\Gamma\longrightarrow\pslr2,
\end{equation}
which satisfy
\begin{equation}\label{eq:equivariance}
\textnormal{dev}(\gamma x)=\rho(\gamma)\,\textnormal{dev}(x) \qquad \textnormal{for every }\gamma\in\Gamma\textnormal{ and }x\in\widetilde S.
\end{equation}

As already alluded above, a branched hyperbolic structure on \(S\) naturally comes with an underlying complex structure \(X\). Such structure lifts to the universal cover, thereby providing a conformal identification \(\widetilde S\cong \mathbb H\), and the developing map \(\textnormal{dev}\) is holomorphic with respect to the lifted complex structure. At a lift of a point of branching order \(k\), it has local degree \(k+1\). Thus it is locally a branched covering, although it need not be a covering of the whole hyperbolic plane.

To construct the developing--holonomy pair, one may first analytically continue the hyperbolic charts on the complement of the branch set. This gives a developing map and a holonomy representation defined a priori on the fundamental group of the punctured surface. Near a branch point of order \(k\), however, the local developing map is of the form \(z\mapsto z^{k+1}\), and therefore has trivial holonomy around a small peripheral loop. Consequently, the holonomy representation factors through \(\pi_1(S)\), while the developing map extends across every lift of a branch point. We also notice that a different choice of the initial chart replaces the pair by
\begin{equation}\label{eq:conjugate-pair}
    (\textnormal{dev},\rho) \longmapsto\bigl(g\circ\textnormal{dev},\,\textnormal{Ad}_g(\,\rho\,) \bigr),\qquad g\in\pslr2.
\end{equation}

As a consequence, the holonomy is intrinsically defined up to conjugation. We write \(\hol(\sigma)=[\rho]\). Whenever we identify a holonomy representation with a prescribed homomorphism, this is understood up to this conjugation, which can be removed by changing the initial developing chart.

Conversely, suppose that \(S\) carries a complex structure and that a non-constant holomorphic map, say \(\textnormal{dev}\colon\widetilde S\longrightarrow\HH^2\) is equivariant with respect to \(\rho\). Its critical set is \(\Gamma\)-invariant and is locally the zero set of a non-zero holomorphic derivative. It therefore descends to a closed discrete subset of the compact surface \(S\), and this subset is finite. The pull-back \(\textnormal{dev}^*h_{\HH}\) descends to a branched hyperbolic structure on \(S\). Its branching orders are the vanishing orders of the differential of \(\textnormal{dev}\), and its holonomy is \([\rho]\).

\begin{defn}\label{def: geo rep}
    A representation \(\rho\colon\Gamma\to\pslr2\) is \emph{geometrisable} if it is the holonomy of a branched hyperbolic structure compatible with the given orientation of \(S\). Equivalently, a representation is geometrisable if there exists a non-constant holomorphic \(\rho\)-equivariant map \(u\colon\widetilde S\longrightarrow\HH^2\).
\end{defn}


\subsection{Moduli space of branched hyperbolic structures and the holonomy map}\label{subsec:moduli-bhs} We next recall the notion of deformation space of branched hyperbolic structures and make explicit the role of the marking. If \((X,\sigma)\) and \((X',\sigma')\) are branched hyperbolic surfaces, an \emph{isomorphism} between them will mean an orientation-preserving isometry \(h\colon(X,\sigma)\to(X',\sigma')\). Equivalently, in branched hyperbolic charts, \(h\) is locally induced by an element of \(\pslr2\). In what follows, we shall denote by \(\operatorname{Isom}^+(\sigma)\) the group of orientation-preserving isometries of \((X,\sigma)\).

\begin{defn}
    A \emph{marked branched hyperbolic structure} on \(S\) is a triple \((X,\sigma,f)\), where \(X\) is a closed oriented surface, \(\sigma\) is a branched hyperbolic structure on \(X\), and \(f\colon S\to X\) is an orientation-preserving diffeomorphism. Two marked structures \((X_1,\sigma_1,f_1)\) and \((X_2,\sigma_2,f_2)\) are equivalent if there exists an orientation-preserving isometry, say \(h\colon(X_1,\sigma_1)\longrightarrow(X_2,\sigma_2)\) such that \(h\circ f_1\) is isotopic to \(f_2\) through orientation-preserving diffeomorphisms. We denote by \(\mathcal{BH}(S)\) the space of marked branched hyperbolic structures on \(S\) up to this equivalence relation.
\end{defn}

\subsubsection{Branched hyperolic structures as developing--holonomy pairs} By using the marking, every branched hyperbolic structure may be pulled back to the fixed topological surface \(S\). Thus \(\mathcal{BH}(S)\) can equivalently be regarded as the space of branched hyperbolic structures on \(S\) modulo the action of \(\operatorname{Diff}_0^+(S)\), the group of orientation-preserving diffeomorphisms isotopic to the identity. As already discussed in \S\ref{subsec:developing-holonomy}, every marked branched hyperbolic structure determines an equivalence class of developing--holonomy pairs. Conversely, such a pair determines the corresponding branched hyperbolic structure. Thus the deformation space admits the following equivalent description:
\begin{equation}\label{eq:identification-bhs}
    \mathcal{BH}(S)
    =
    \left\{\,\,\,
    \begin{gathered}
        \text{branched hyperbolic}\\
        \text{structures on } S
    \end{gathered}\,\,\,
    \right\}
    \longleftrightarrow
    \Big\{
        \big[\,g\circ\textnormal{dev},\,\operatorname{Ad}_g(\,\rho\,)\,\big]
        \,\big|\,
        g\in\pslr2
    \Big\}_{\sim\,\operatorname{Diff}_0^+(S)}.
\end{equation}

On the other hand, changing the marking by a diffeomorphism isotopic to the identity corresponds to precomposing the developing map by a lift of that diffeomorphism. Consequently, the deformation space \(\mathcal{BH}(S)\) may equivalently be described as the corresponding quotient of developing--holonomy pairs by these two actions. The marking will be kept throughout the paper. In particular, constructing a branched hyperbolic structure whose holonomy has the same image as a prescribed representation \(\rho\) is not sufficient to realise \(\rho\). In fact, what must agree, up to conjugation in \(\pslr2\), is the homomorphism induced from the marked surface group \(\pi_1(S)\). This distinction will be important when keeping track of the prescribed representation through the pinching and covering constructions used below.

\subsubsection{Representation space and holonomy map}\label{sssec:holomap} We first recall some basic facts about surface group representations into \(\pslr2\). Let \(S\) be a closed oriented surface of genus \(g\geq2\). In our setting, the \emph{representation space} is \(\operatorname{Hom}\bigl(\pi_1(S),\pslr2\bigr)\) endowed with its natural structure of a real analytic variety. The group \(\pslr2\) acts on this space by conjugation, and the corresponding \emph{character variety} is the quotient
\begin{equation}\label{eq:character-variety}
    \mathfrak X(S)
    \defeq
    \operatorname{Hom}\bigl(\pi_1(S),\pslr2\bigr)\big/\pslr2.
\end{equation}

In his seminal work \cite{Goldman88}, Goldman showed that both the representation space and the associated character variety have \(4g-3\) connected components, distinguished by the Euler number, a topological invariant naturally associated with a representation. For our convenience, we recall the definition and the main properties of the Euler number in the following subsection; see \S\ref{subsec:euler-number}. Finally, we observe that by equipping the space of developing maps with the compact-open topology and the representation space with its natural topology gives \(\mathcal{BH}(S)\) the corresponding quotient topology. The holonomy construction then defines a continuous map, called \textit{holonomy map},
\begin{equation}\label{eq:holonomy-map}
    \hol\colon
    \mathcal{BH}(S)
    \longrightarrow
    \mathfrak X(S),
\end{equation}
which associates to a marked branched hyperbolic structure the conjugacy class of its holonomy representation. Under this perspective, our main result characterises which discrete representations arise as the holonomy of some branched hyperbolic structure.

\medskip

\subsection{Euler number}\label{subsec:euler-number} The action of \(\pslr2\) on the oriented circle at infinity \(\partial\HH^2\cong\mathbb{RP}^1\) associates to any representation \(\rho\) a flat oriented circle bundle
\begin{equation}\label{eq:flat-circle-bundle}
    E_\rho=\frac{\widetilde S\times\partial\HH^2}{\Gamma}\longrightarrow S, \qquad \gamma\cdot(x,\xi)=\bigl(\gamma x,\rho(\gamma)\xi\bigr).
\end{equation}

Its Euler class belongs to \(H^2(S;\Z)\cong\Z\) and it can be seen as the obstruction to the existence of global sections of \(E_\rho\). We define the \emph{Euler number} of \(\rho\) by
\begin{equation}\label{eq:euler-number}
    e(\rho)=\bigl\langle e(E_\rho),[S]\bigr\rangle\in\Z.
\end{equation}

\textit{Convention:} We orient \(\partial\HH^2\) so that the endpoint map from the positively oriented circle of unit tangent vectors at any point of \(\HH^2\) is orientation preserving.

We next give an equivalent description of the Euler class in terms of tangent vectors. This will also provide the circle bundle used later in the proof of the main theorem. For this purpose, choose a smooth \(\rho\)-equivariant map \(u\colon\widetilde S\to\HH^2\). Such a map exists because the associated bundle with fibre \(\HH^2\) has contractible fibres. Define an oriented real plane bundle and its unit-circle bundle by
\begin{equation}
\label{eq:associated-plane-bundle}
V_u=\bigl(u^*T\HH^2\bigr)/\Gamma,
\qquad
M_u=\bigl(u^*T^1\HH^2\bigr)/\Gamma.
\end{equation}
The action on the tangent vectors is given by the differential of \(\rho(\gamma)\). By associating a unit vector to the endpoint of the geodesic ray that it determines, the resulting map provides an orientation-preserving bundle isomorphism \(M_u\cong E_\rho\). In particular, \(e(V_u)=e(E_\rho)\). We may notice that the bundle \(M_u\) is a principal circle bundle under rotation of the unit tangent vectors. The flat structure on \(E_\rho\), however, is a flat structure with group \(\pslr2\); it is not being identified with a flat principal-circle connection. Our sign convention can now be read off from a hyperbolic surface. If \(j\colon\Gamma\to\pslr2\) is the holonomy of a hyperbolic metric compatible with the orientation of \(S\), the differential of its developing map identifies the corresponding bundle \(V_{\textnormal{dev}}\) with \(TS\). Hence
\begin{equation} \label{eq:fuchsian-sign}
e(\,j\,)=e(\,TS\,)=\chi(S)=2-2g<0.
\end{equation}

With this convention the Milnor--Wood inequality reads \(\lvert e(\rho)\rvert\leq 2g-2\); see \cite{Milnor1958,Wood1971}. Naturality of the Euler class will be used repeatedly. If \(f\colon S'\to S\) is a map of closed oriented surfaces of degree \(n\), then
\begin{equation}
\label{eq:euler-naturality}
e(\rho\circ f_*)=n\,e(\rho).
\end{equation}
Indeed, \(E_{\rho\circ f_*}\) is isomorphic to \(f^*E_\rho\), and evaluation on the fundamental class gives \eqref{eq:euler-naturality}. In particular, restriction to a subgroup corresponding to an oriented covering of degree \(n\) multiplies the Euler number by \(n\).

\begin{lem}\label{lem:euler-branching}
    If \(\rho\) is the holonomy of a branched hyperbolic structure \(\sigma\), then
    \begin{equation}\label{eq:euler-branching}
        e(\rho)=\chi(S)+\kappa(\sigma).
    \end{equation}
\end{lem}

For the readers' convenience, we provide a concise proof of Lemma \ref{lem:euler-branching}. Alternative arguments with more geometric-topological flavour may be found in \cite{Faraco2021} and \cite{Mathews}.

\begin{proof}
    Since \(\rho\) is the holonomy of the branched hyperbolic structure \(\sigma\), there exists a developing map \(\textnormal{dev}\colon\widetilde S\longrightarrow\HH^2\) which is \(\rho\)-equivariant by Definition \ref{def: geo rep}. We may therefore take \(u=\textnormal{dev}\) in the construction of \eqref{eq:associated-plane-bundle}. Since \(\sigma\) determines a complex structure on \(S\) for which \(\textnormal{dev}\) is holomorphic, the bundle \(V_{\textnormal{dev}}\) inherits the structure of a holomorphic line bundle over \(S\). The differential of the developing map descends to a holomorphic section of \(\operatorname{Hom}_{\C}(TS,V_{\textnormal{dev}})\). Near a branch point of order \(k\), this section vanishes to order \(k\), and hence its zero divisor is precisely \(D_\sigma\). Taking first Chern numbers gives
    \begin{equation}
        c_1(V_{\textnormal{dev}})[S]-c_1(TS)[S]=\deg D_\sigma.
    \end{equation}
    Since these two Chern numbers are respectively \(e(\rho)\) and \(\chi(S)\), the claimed identity follows.
\end{proof}

\smallskip

\subsection{Volume}\label{subsec:volume} We next recall the volume of a representation. This invariant provides the natural link between the topology of a representation and the hyperbolic geometry of an equivariant map. Let \(\rho\colon\Gamma\to\pslr2\) be a representation. Since \(\HH^2\) is contractible, there exists a \(\rho\)-equivariant smooth map \(u\colon\widetilde S\longrightarrow\HH^2\). Let 
\[
    \omega_{\HH}=\frac{dx\wedge dy}{y^2}
\]
denote the oriented hyperbolic area form on \(\HH^2\). The \(\rho\)-equivariance of the mapping \(u\), together with the invariance of \(\omega_{\HH}\) under \(\pslr2\), implies that \(u^*\omega_{\HH}\) is \(\Gamma\)-invariant and therefore descends to a \(2\)-form on \(S\). Its integral over \(S\) will define the volume of \(\rho\).

\begin{defn}
Let \(\rho\colon\Gamma\to\pslr2\) be a representation. The \emph{volume} of \(\rho\) is defined as the value of the integral
\begin{equation}\label{eq:volume-definition}
    \textnormal{vol}(\rho)=\int_S u^*\omega_{\HH},
\end{equation}
where the integrand denotes the descended form on \(S\) and the integral uses the prescribed orientation of \(S\).
\end{defn}

\begin{rmk}
    The value in \eqref{eq:volume-definition} does not depend on \(u\). In fact, two equivariant maps are equivariantly homotopic by pointwise geodesic interpolation in \(\HH^2\). The pull-back of \(\omega_{\HH}\) under this homotopy descends to a closed two-form on \(S\times[0,1]\), and Stokes' theorem shows that the integrals at the two ends agree.
\end{rmk}

The relation between volume and Euler number is explicit and, according to our convention, given by the following identity. 

\begin{lem}\label{lem:volume-euler}
    For a representation \(\rho\colon\Gamma\longrightarrow\pslr2\) with volume \(\textnormal{vol}(\rho)\) and Euler number \(e(\rho)\) the following identity holds:
    \begin{equation}
        \textnormal{vol}(\,\rho\,)=-2\pi e(\,\rho\,).
    \end{equation}
\end{lem}

For completeness, we provide a proof for the readers' convenience.

\begin{proof}
    Let \(\nabla^{\HH}\) denote the Levi--Civita connection of \(T\HH^2\), and let \(u^*\nabla^{\HH}\) denote its pull-back connection on \(u^*T\HH^2\) via the \(\rho\)-equivariant map \(u\colon\widetilde S\longrightarrow\HH^2\). Since every element of \(\pslr2\) acts on \(\HH^2\) by an orientation-preserving isometry, it preserves the Levi--Civita connection. Since the map \(u\) is \(\rho\)-equivariant, the pull-back connection \(u^*\nabla^{\HH}\) is preserved by the natural \(\Gamma\)-action on \(u^*T\HH^2\). Hence \(u^*\nabla^{\HH}\) descends to a connection \(\nabla^u\) on the quotient bundle \(V_u\to S\). We now compute its Euler form. For an oriented Riemannian surface, the Euler form of the Levi--Civita connection is given by
    \begin{equation}
        \frac{K}{2\pi}\,dA,
    \end{equation}
    where \(K\) is the Gaussian curvature and \(dA\) is the oriented area form. Since the Gaussian curvature of the hyperbolic plane is constant and equal to \(K=-1\), the Euler form of \(T\HH^2\) is
    \begin{equation}
        -\frac{1}{2\pi}\omega_{\HH}.
    \end{equation}
    By naturality under pull-back, the Euler form of the descended connection \(\nabla^u\) on \(V_u\) is therefore represented by
    \begin{equation}
        -\frac{1}{2\pi}u^*\omega_{\HH}.
    \end{equation}
    
    By Chern--Weil theory, the integral of this Euler form over \(S\) is equal to the Euler number of the bundle \(V_u\). Since \(e(V_u)=e(E_\rho)\), and hence \(e(V_u)[S]=e(\rho)\), we readily obtain
    \begin{equation}
        e(\rho)=\int_S-\frac{1}{2\pi}u^*\omega_{\HH}=-\frac{1}{2\pi}\textnormal{vol}(\rho).
    \end{equation}
    
    Therefore identity \( \textnormal{vol}(\rho)=-2\pi e(\rho)\) follows as desired.
\end{proof}


If \([\rho]=\hol(\sigma)\), one may use the developing map to compute the volume. The pull-back area form is the area form of \(\sigma\) away from the branch points, and the branch points have measure zero. Thus
\begin{equation}
\label{eq:volume-geometric}
\textnormal{vol}(\,\rho\,)
=\textnormal{Area}(S,\sigma)
=2\pi\bigl(2g-2-\kappa(\sigma)\bigr)>0.
\end{equation}
In particular, positive volume is necessary for geometrisability with the chosen orientation. Equations \eqref{eq:euler-naturality} and Lemma \ref{lem:volume-euler} also give
\begin{equation}
\label{eq:volume-naturality}
\textnormal{vol}(\,\rho\circ a_*\,)
=\deg(\,a\,)\,\textnormal{vol}(\,\rho\,)
\end{equation}
for maps \(a\) of closed oriented surfaces.

\smallskip

\section{Seifert fibrations and fibre-preserving maps}\label{sec:seifert-fibrations}

The proof of the main theorem will make essential use of maps between Seifert fibred \(3\)-manifolds. We therefore collect in this section the basic facts that will be needed later. In particular, we distinguish principal circle bundles from Seifert fibrations with exceptional fibres, recall the local model around an exceptional fibre, and describe the behaviour of fibre-preserving branched coverings in such neighbourhoods. We also record the corresponding information on fundamental groups.

\begin{rmk}
    We stress that the Seifert fibred \(3\)-manifolds considered below should not, in general, be regarded as manifolds modelled on \(\HH^2\times\R\). For Seifert manifolds with hyperbolic base orbifold, the relevant Thurston geometry depends on the Euler class of the fibration. In particular, the unit tangent bundle of a hyperbolic surface or orbifold naturally carries a \(\widetilde{\pslr2}\)-geometry, whereas the geometry \(\HH^2\times\R\) occurs when the corresponding Seifert fibration has vanishing Euler class. In what follows, however, we shall only use the topological and fibred structure of these manifolds.
\end{rmk}

\subsection{Regular and exceptional fibres}

A Seifert fibration of an oriented \(3\)-manifold \(M\) is a decomposition of \(M\) into circles, called \emph{fibres}, such that every fibre admits a fibred solid torus as a neighbourhood. All but finitely many fibres are regular; the remaining ones are called \emph{exceptional}. The space of fibres is naturally an oriented \(2\)-orbifold: regular fibres project to regular points, whereas an exceptional fibre of multiplicity \(m\) projects to a cone point of order \(m\).

A principal circle bundle over a smooth surface is the special case in which every fibre is regular. In particular, over a disc \(D\) it is trivial, \(D\times\sph^1\to D\). An exceptional fibre has a different local model. Let \(m\geq2\), let \(q\ge1\) such that \(\gcd(q,m)=1\), and set
\begin{equation}
    \zeta=e^{2\pi i/m}.
\end{equation}
The cyclic group \(\Z_m\) acts freely on \(D\times\sph^1\) by
\begin{equation}
    (z,u)\longmapsto(\zeta z,\zeta^q u).
\end{equation}
The quotient
\begin{equation}\label{eq:exceptional-solid-torus}
    V_{m,q}=\frac{D\times\sph^1}{\Z_m}
\end{equation}
is a smooth solid torus carrying a Seifert fibration over the orbifold disc \(D/\Z_m\). The image of the circle \(\{0\}\times\sph^1\) is an exceptional fibre of multiplicity \(m\).

Equivalently, one may cut \(V_{m,q}\) along a transverse disc and obtain \(D\times[0,1]\), whose two ends are reglued after a rotation of order \(m\). The central fibre closes after one passage, whereas a nearby regular fibre closes only after \(m\) passages. The covering
\begin{equation}\label{eq:unwrapping-exceptional-fibre}
    D\times\sph^1\longrightarrow V_{m,q}
\end{equation}
has degree \(m\) and unwraps the exceptional fibre.

\subsection{The fibre slope}\label{ssec:slope} The multiplicity of a fibre can also be read on the boundary torus. Let \(V\) be a fibred solid torus, let \(\mu\) be a meridian, and choose a longitude \(\lambda\) so that \(H_1(\partial V;\Z)=\Z\mu\oplus\Z\lambda\). If \(h\) is the homology class of a regular fibre on \(\partial V\), then
\begin{equation}\label{eq:seifert-fibre-slope}
    h=a\mu+m\lambda,
    \qquad
    \gcd(a,m)=1.
\end{equation}
Hence the multiplicity of the core fibre is \(m=\lvert\mu\cdot h\rvert\), where \(\mu\cdot h\) denotes the algebraic intersection number on \(\partial V\). In particular, the core fibre is regular precisely when \(m=1\). We shall use this characterisation when studying vertical pinches; compare \cite[Lemma~3.2.8]{Konstantinov2019}.

\subsection{Fibre-preserving maps near an exceptional fibre}\label{sssec:fibre-preserving-maps} Let \(B\colon V_1\to V_2\) be a fibre-preserving branched covering. Assume that the core fibre of \(V_1\) is regular, that the core fibre of \(V_2\) is exceptional of multiplicity \(m\), and that \(B\) has degree \(+1\) on regular fibres. Thus the source is locally a fibre product \(V_1\cong D_w\times\sph^1\), where \(D_w\) is a disc with coordinate \(w\), whereas the target is modelled on
\begin{equation}
    V_2=V_{m,q}=(D_z\times\sph^1)/\Z_m.
\end{equation}

Let
\begin{equation}\label{eq:unwrapping-map}
    p\colon\widehat V_2\defeq D_z\times\sph^1\longrightarrow V_2
\end{equation}
be the \(m\)-fold covering introduced in \eqref{eq:unwrapping-exceptional-fibre}. The Seifert fibration of \(\widehat V_2\) is a product fibration, so the exceptional fibre of \(V_2\) becomes a regular fibre after passing to this covering. Let \(h_1\) and \(h_2\) denote the classes of regular fibres of \(V_1\) and \(V_2\), respectively. Since \(B\) has regular-fibre degree \(+1\),
\begin{equation}
    B_*(h_1)=h_2.
\end{equation}
On the other hand, the subgroup \(p_*(\pi_1(\widehat V_2))\subset\pi_1(V_2)\) is precisely the subgroup generated by \(h_2\). Hence the covering-space lifting criterion gives a lift, say \( \widehat B\colon V_1\longrightarrow\widehat V_2\), such that
\begin{equation}
    B=p\circ\widehat B.
\end{equation}

Both \(V_1\) and \(\widehat V_2\) are now product fibrations. Denote their respective projections by \(q_1\colon V_1\longrightarrow D_w\) and \(\widehat q_2\colon\widehat V_2\longrightarrow D_z\). Since \(\widehat B\) is fibre preserving and has degree \(+1\) on the fibres, it induces a branched map, say \(\widehat b\colon D_w\longrightarrow D_z\). After choosing oriented coordinates centred at the two core fibres, this map has the local form
\begin{equation}\label{eq:transverse-ramification}
    z=w^r
\end{equation}
for some \(r\geq1\). We call \(r\) the \emph{transverse ramification index} of \(B\). Equivalently, \(r\) is the local degree of the mapping \(B\) between sufficiently small discs transverse to the source and target core fibres. Thus \(r\) measures the genuine branching of the three-dimensional map in the direction transverse to the fibres. The relation between these maps and the corresponding base spaces is summarised by the commutative diagram
\begin{equation}\label{eq:local-seifert-diagram}
    \begin{tikzcd}[column sep=large,row sep=large]
        V_1
        \arrow[r,"\widehat B"]
        &
        \widehat V_2
        \arrow[r,"p"]
        &
        V_2
        \\
        D_w
        \arrow[r,"\widehat b"']
        &
        D_z
        \arrow[r,"\pi_m"']
        &
        \lvert D_z/\Z_m\rvert
        \arrow[from=1-1,to=2-1,"q_1"']
        \arrow[from=1-2,to=2-2,"\widehat q_2"']
        \arrow[from=1-3,to=2-3,"q_2"]
    \end{tikzcd}
\end{equation}

where \(q_2\) is the Seifert projection and \(\pi_m\) is the quotient map on the base. If \(u\) denotes a coordinate on the underlying topological disc \(|D_z/\Z_m|\), then
\begin{equation}\label{eq:orbifold-quotient-coordinate}
    u=z^m,
\end{equation}
and \eqref{eq:transverse-ramification} and \eqref{eq:orbifold-quotient-coordinate} combined together give
\begin{equation}
    u=w^{mr}.
\end{equation}
Consequently, although the transverse ramification index of \(B\) is \(r\), the local degree of the induced map on the underlying base surfaces is \(mr\). If the target base carries its natural orbifold cone angle \(2\pi/m\), its pull-back angle is therefore
\begin{equation}\label{eq:angle-2pir-preliminaries}
    mr\cdot\frac{2\pi}{m}=2\pi r.
\end{equation}
Thus the factor \(r\) measures the actual transverse branching, whereas the additional factor \(m\) comes from the exceptional multiplicity of the target fibre.

\begin{rmk}\label{rmk:transverse-disc-exceptional}
    For the readers' convenience, we stress a quite subtle distinction between the covering \(p\colon\widehat V_2\to V_2\) and the quotient map \(\pi_m\colon D_z\to |D_z/\Z_m|\). The former is an ordinary unbranched covering of smooth \(3\)-manifolds. In particular, if \(\widehat\Delta\subset\widehat V_2\) is a sufficiently small disc transverse to the central fibre, then the restriction \(p|_{\widehat\Delta}\) is an embedding onto a transverse disc \(\Delta\subset V_2\). Thus the factor \(m\) should not be interpreted as ramification of \(p\) in the transverse direction. What is different near an exceptional fibre is that a transverse disc \(\Delta\subset V_2\) is not a local section of the Seifert projection. Its projection to the underlying base has local degree \(m\). More precisely, after identifying \(\widehat\Delta\) with \(D_z\), one has
    \begin{equation}
        q_2\circ p|_{\widehat\Delta}=\pi_m,
    \end{equation}
    and \(\pi_m\) is locally given by \(u=z^m\). Hence, when \(r=1\), the lifted map \(\widehat B\) is locally unramified transversely and \(B\) maps a transverse disc locally homeomorphically onto a transverse disc of \(V_2\); nevertheless, the induced map on the underlying base surfaces has local degree \(m\). The factor \(m\) is therefore entirely due to the exceptional Seifert structure of the target.
\end{rmk}

\smallskip

\subsection{Fundamental groups} Finally, let \(\sph^1\to M\to S\) be an oriented principal circle bundle over an aspherical closed oriented surface. The homotopy exact sequence gives
\begin{equation}\label{eq:circle-bundle-central-extension}
    1\longrightarrow\langle \,h\,\rangle
    \longrightarrow\pi_1(M)
    \longrightarrow\pi_1(S)
    \longrightarrow1,
    \qquad
    \langle \,h\,\rangle\cong\Z,
\end{equation}
where \(h\) denotes the positively oriented fibre class. The fibre subgroup is central. For an oriented Seifert fibration \(M\to\mathcal O\) over an orientable hyperbolic orbifold \(\mathcal O\), the regular fibre similarly determines a central infinite cyclic subgroup and
\begin{equation}\label{eq:seifert-orbifold-extension}
    1\longrightarrow\langle \,h\,\rangle
    \longrightarrow\pi_1(M)
    \longrightarrow\pi_1^{\mathrm{orb}}(\mathcal O)
    \longrightarrow1.
\end{equation}
These two exact sequences will be used repeatedly below.

\smallskip

\section{Geometrisation of discrete representations}\label{sec:discrete-geometrisation}

This section is devoted to the proof of the main result, \textit{i.e.}, Theorem \ref{thm:discrete-geometrisation}. In what follows we shall say that a representation is said to be \emph{purely hyperbolic} if every non-identity element of its image is hyperbolic. The geometrisation problem for these representations was treated in \cite[Theorem 3.14]{Faraco2021}. We shall use that result for the torsion-free case. Representations with discrete image containing a parabolic element have Euler number zero and hence do not occur under our positive-volume hypothesis. We first record the stronger statement needed to make this reduction complete; compare \cite[Theorem 1.2]{FunarWolff2007}.

\begin{lem}
\label{lem:noncompact-euler-zero}
Let \(\rho\colon\Gamma\to\pslr2\) have discrete image \(G\). If \(G\backslash\HH^2\) is non-compact, then \(e(\rho)=0\).
\end{lem}

\begin{proof}
    Since the group \(G\) is finitely generated then there exists a torsion-free subgroup \(H\) of finite index \(n\) by Selberg's Lemma. Set \(K=\rho^{-1}(H)\), and let \(S_H\longrightarrow S\) be the covering map corresponding to \(K\). Notice that, since \(\rho\) is surjective onto \(G\), this covering has degree \(n\). The quotient space \(\Sigma=H\backslash\HH^2\) is a smooth non-compact aspherical surface, being a finite covering of the non-compact orbifold \(G\backslash\HH^2\). In particular, \(\Sigma\) is a \(K(H,1)\) and \(H^2(\Sigma;\Z)=0\). Hence the Euler class of the inclusion \(H\hookrightarrow\pslr2\), viewed as a class in \(H^2(\Sigma;\Z)\), vanishes. The restricted representation \(\rho_H\colon K\to H\hookrightarrow\pslr2\) pulls back this class to \(S_H\), and therefore \(e(\rho_H)=0\). On the other hand, by \eqref{eq:euler-naturality}, 
    \begin{equation}
        n\,e(\rho)=e(\rho_H)=0.
    \end{equation}
    Since \(n>0\), this proves the assertion. The same argument includes elementary images; in particular, \(H\) is allowed to be trivial.
\end{proof}

\begin{rmk}\label{rmk:elementary-euler-zero}
    Every representation \(\rho\colon\pi_1(S)\to\pslr2\) with discrete elementary image has vanishing Euler number. Indeed, a discrete elementary subgroup of \(\pslr2\) is virtually cyclic. Hence there exists a finite-index subgroup \(H<\rho(\pi_1(S))\) which is cyclic. Let \(S_H\to S\) be the finite covering corresponding to \(\rho^{-1}(H)\), and let \(\rho_H\) denote the restriction of \(\rho\). If \(H\) is infinite cyclic, then \(\rho_H\) factors through \(\Z\), and therefore its Euler class vanishes since \(H^2(\Z;\Z)=0\). If \(H\) is finite, we may pass to a further finite cover on which the representation is trivial. In either case, the Euler number vanishes on some finite cover. By naturality of the Euler class,
    \begin{equation}
        e(\rho_H)=[\,\rho(\pi_1(S)):H\,]\,e(\rho),
    \end{equation}
    and hence \(e(\rho)=0\). In particular, elementary representations have null volume.
\end{rmk}

It follows that a discrete representation of positive volume has cocompact image. Such an image contains no parabolic elements. If it is torsion-free, it is purely hyperbolic. Notice that absence of parabolic elements alone would not imply cocompactness; the volume assumption is used through Lemma \ref{lem:noncompact-euler-zero}. We are therefore left to prove the following assertion.

\begin{prop}\label{prop:torsion-geometrisation}
    Let \(\rho\colon\Gamma\longrightarrow\pslr2\) have discrete cocompact image \(G\) containing torsion and assume \(\textnormal{vol}(\rho)>0\). Then \(\rho\) is geometrisable.
\end{prop}

We shall prove Proposition \ref{prop:torsion-geometrisation} into several steps collected in \S\S\ref{subsec:circle-bundles}-\ref{subsec:completion}. The main construction replaces the orbifold target by its unit tangent bundle and applies Rong's factorisation theorem, see \cite{Rong1993}, to a map of Seifert fibred three-manifolds. We then recover a branched hyperbolic structure on a surface and reinsert the handles removed by the factorisation. Throughout the construction, the image group is the actual subgroup \(G=\rho(\Gamma)\), and the resulting homomorphisms retain the marking of the source.

\subsection{The associated circle bundles}\label{subsec:circle-bundles} Let \(\rho\colon\Gamma\to\pslr2\) be a discrete representation with positive volume, \textit{i.e.}, \(\textnormal{vol}(\rho)>0\). Assume \(G\) is a cocompact group with torsion. Let \(u\colon\widetilde S\longrightarrow\HH^2\) be a smooth \(\rho\)-equivariant map, as in \S\ref{subsec:euler-number}, and consider
\begin{equation}\label{eq:three-manifolds}
    M=\bigl(u^*T^1\HH^2\bigr)/\Gamma, \qquad N=G\backslash T^1\HH^2.
\end{equation}

Notice that first projection \(q_M\colon M\to S\) determines an oriented principal circle bundle over \(S\). In the definition of \(M\), a point before taking the quotient is a pair \((x,v)\), where \(v\) is a unit tangent vector based at \(u(x)\), and
\begin{equation}\label{eq:action-on-pullback}
    \gamma\cdot(x,v)=\bigl(\gamma x,D\rho(\gamma)_{u(x)}v\bigr).
\end{equation}
We may observe that the fibre-wise rotation action of \(\sph1\) commutes with the \(\Gamma\)-action, since every element of \(\pslr2\) preserves both the hyperbolic metric and the orientation. Hence the \(\sph1\)-action descends to \(M\), where it remains free and transitive on the fibres of \(q_M\colon M\longrightarrow S\). Therefore \(q_M\) is a principal circle bundle.

\begin{rmk}\label{rmk:source-target-fibrations}
    There is an important distinction between the two Seifert fibrations introduced above. The manifold \(M\to S\) is a principal circle bundle, and therefore \textit{all} its fibres are regular. On the other hand, \(N\to\mathcal O\) is in general only a Seifert fibration: the exceptional fibres of \(N\) correspond precisely to the cone points of the orbifold \(\mathcal O\). This distinction will be important when applying Rong's factorisation in \S\ref{subsec:rong}.
\end{rmk}

Let \(\mathcal O=G\backslash\HH^2\). Although \(\mathcal O\) may have cone points, \(N\) is a smooth three-manifold. Indeed, an orientation-preserving hyperbolic isometry that fixes a unit tangent vector fixes its base point and has trivial derivative there, and is consequently the identity. Thus the discrete group \(G\) acts freely and properly discontinuously on \(T^1\HH^2\). The projection \(q_N\colon N\longrightarrow\mathcal O\) is a Seifert fibration and a cone point of order \(m\) in \(\mathcal O\) corresponds to an exceptional fibre of multiplicity \(m\).

The manifolds \(M\) and \(N\) are closed, connected and orientable. For later purposes, we observe that they are both aspherical and have infinite fundamental groups, as required for the application of Rong's factorisation result in \S\ref{subsec:rong}. Indeed, pulling the circle bundle \(M\to S\) back to the contractible surface \(\widetilde S\) gives a trivial bundle \(\widetilde S\times\sph^1\), which covers \(M\). The universal cover of \(M\) is therefore homeomorphic to \(\widetilde S\times\R\). In the same fashion, since the action of \(G\) on \(T^1\HH^2\) is free, the projection map \(T^1\HH^2\to N\) is a covering. Since \(T^1\HH^2\cong\HH^2\times\sph^1\), it readily follows that the universal cover of \(N\) is homeomorphic to \(\HH^2\times\R\). Both universal covers are contractible, thereby proving the asphericity. Moreover, they are non-compact, so the compactness of \(M\) and \(N\) implies that their fundamental groups are infinite. 

There is a natural map \(\widetilde F\colon u^*T^1\HH^2\to T^1\HH^2\) that associated with every pair \((x,v)\) the unit tangent vector \(v\); notice that the base point is uniquely determined by \(v\). By \eqref{eq:action-on-pullback}, the action \(\Gamma\) on the pull-back bundle satisfies \(\widetilde F\big(\,\gamma\cdot(x,v)\,\big)=D\rho(\,\gamma\,)_{u(x)}v\). As \(\rho(\,\gamma\,)\in G\), the vectors \(v\) and \(D\rho(\,\gamma\,)_{u(x)}v\) belong to the same \(G\)-orbit in \(T^1\HH^2\). Consequently, \(\widetilde F\) descends to a well-defined map 
\begin{equation} \label{eq:map-F}
    F\colon M\longrightarrow N, \qquad F([x,v])=[\,v\,].
\end{equation}

The fibrewise rotation actions descend to both \(M\) and \(N\). By the definition of \(F\), for every \(\theta\in\sph^1\) the identity
\begin{equation}
    F\bigl([x,R_\theta v]\bigr)
    =
    [R_\theta v]
    =
    R_\theta\cdot F([x,v])
\end{equation}
holds. Thus \(F\) is \(\sph^1\)-equivariant and, in particular, maps fibres of \(q_M\) to fibres of \(q_N\). If the image fibre is regular, its stabiliser under the rotation action is trivial. Since the rotation action on the source fibre is also free, the restriction of \(F\) to that fibre is an orientation-preserving map of degree one. Suppose instead that the image fibre is exceptional of multiplicity \(m\). Its stabiliser under the rotation action is then the cyclic subgroup \(\Z_m<\sph^1\). Consequently, the restriction of \(F\) to the source fibre is identified with the quotient map \(\sph^1\longrightarrow \sph^1/\Z_m\). After identifying \(\sph^1/\Z_m\) with the underlying circle of the exceptional fibre, this map has degree \(m\).


\begin{rmk}
    The degree \(m\) of the restriction of \(F\) to a source fibre mapping onto an exceptional fibre can also be read from the fibre-slope description in \S\ref{ssec:slope}. Let \(V\) be a fibred solid torus neighbourhood of the exceptional fibre, and let \(c\in H_1(V;\Z)\) be the primitive class represented by its core. Let \(\imath\colon \partial V\hookrightarrow V\) be the natural inclusion. Every regular fibre contained in \(\partial V\) represents the same slope \(h\in H_1(\partial V;\Z)\), which has the form \(h=a\mu+m\lambda\). Here \(\mu\) is a meridian and \(\lambda\) is oriented so that its image under the inclusion \(\imath_*\colon H_1(\partial V;\Z)\to H_1(V;\Z)\cong H_1(\partial V;\Z)/\mathbb Z\mu\) is \(c\). Since the meridian bounds a disc in \(V\), its image is trivial, and hence \(h\) maps to \(mc\). Thus a regular fibre represents the \(m\)-th power of the exceptional core, which is precisely why the restriction of \(F\) to a source fibre has degree \(m\) onto the underlying exceptional fibre.
\end{rmk}


\subsection{Fundamental groups and the degree}\label{subsec:fundamental-groups-degree} Choose base points and paths once and for all. Let \(h_M\) and \(h_N\) denote the homotopy classes of positively oriented regular fibres. The bundle and Seifert fibration give central extensions
\begin{equation}\label{eq:central-extensions}
    \begin{aligned}
        1&\longrightarrow\langle h_M\rangle\longrightarrow\pi_1(M)\longrightarrow\Gamma\longrightarrow1,\\
        1&\longrightarrow\langle h_N\rangle\longrightarrow\pi_1(N)\longrightarrow G\longrightarrow1,
    \end{aligned}
    \qquad
    \langle h_M\rangle\cong\langle h_N\rangle\cong\Z.
\end{equation}
For the second sequence, one may use the regular covering \(T^1\HH^2\longrightarrow N\). Its fundamental group is generated by the positively oriented tangent circle, and its deck group is \(G\). The fibre classes are central because the base and fibre orientations are preserved. The map \(F_*\) sends \(h_M\) to \(h_N\), and the induced map on the quotients in \eqref{eq:central-extensions} is precisely \(\rho\). The latter identification also follows directly from \eqref{eq:action-on-pullback}: lifting a loop to the pull-back of \(T^1\HH^2\) produces the deck transformation \(\rho(\gamma)\).

\begin{lem}\label{lem:F-surjective-degree}
    The homomorphism \(F_*\colon\pi_1(M)\to\pi_1(N)\) is surjective. Moreover,
    \begin{equation}\label{eq:degree-formula}
        \deg F=\frac{e(\,\rho\,)}{\chi_{\mathrm{orb}}(\mathcal O)}=\frac{\textnormal{vol}(\,\rho\,)}{\textnormal{vol}(\mathcal O)}>0.
    \end{equation}
In particular, the degree is a positive integer.
\end{lem}

\begin{proof}
    Let \(\alpha\in\pi_1(\,N\,)\). Since \(\rho\colon\Gamma\to G\) is onto, there is \(\beta\in\pi_1(M)\) whose image \(F_*(\,\beta\,)\) has the same projection to \(G\) as \(\alpha\). Therefore \(\alpha\,F_*(\,\beta\,)^{-1}\) is a power of \(h_N\). The identity \(F_*(h_M)=h_N\) shows that this power also lies in the image of \(F_*\), thereby proving surjectivity. For the degree calculation, choose a torsion-free subgroup \(H<G\) of index \(n\), set \(K=\rho^{-1}(H)\), and write
    \begin{equation}\label{eq:finite-covers}
        p_H\colon S_H\longrightarrow S, \qquad \Sigma=H\backslash\HH^2.
    \end{equation} 

    Both the surface covering \(p_H\) and the orbifold covering \(\Sigma\to\mathcal O\) have degree \(n\). The map \(u\), as introduced in \S\ref{subsec:circle-bundles}, descends to \(u_H\colon S_H\to\Sigma\), inducing the restricted representation \(\rho_H\colon K\to H\). Pulling the circle bundle \(q_M\colon M\longrightarrow S\) back to \(S_H\) via the mapping \(p_H\) yields a circle bundle \(M_H\), and \(F\) lifts to a map \(F_H\) in the commutative diagram
    \begin{equation}\label{eq:degree-covering-diagram}
    \begin{tikzcd}[column sep=large,row sep=large]
        M_H \arrow[r,"F_H"] \arrow[d] & T^1\Sigma \arrow[d]\\
        M \arrow[r,"F"'] & N.
    \end{tikzcd}
    \end{equation}
Each vertical map has degree \(n\). Since \(F_H\) covers \(u_H\) and has oriented fibre degree one, the degree of a map of oriented circle bundles gives
\begin{equation}\label{eq:degree-FH}
    \deg F_H=\deg u_H.
\end{equation}
This equality can also be seen by choosing a regular value of \(u_H\): every point above it contributes its local base degree, multiplied by the positive fibre degree one. Taking degrees in \eqref{eq:degree-covering-diagram} therefore yields \(\deg F=\deg u_H\). The representation \(\rho_H\) thus arises as the composition of \((u_H)_*\) with the orientation-compatible Fuchsian inclusion \(H\hookrightarrow\pslr2\). By \eqref{eq:fuchsian-sign} and \eqref{eq:euler-naturality},
    \begin{equation}\label{eq:euler-degree-cover}
        n\,e(\rho)=e(\rho_H)=\deg(u_H)\,\chi(\Sigma).
    \end{equation}
The orbifold covering satisfies \(\chi(\Sigma)=n\chi_{\mathrm{orb}}(\mathcal O)\). As a consequence, by cancelling \(n\) gives the first equality in \eqref{eq:degree-formula}; whereas Lemma \ref{lem:volume-euler} and \eqref{eq:orbifold-area} combined together yield the second. Both volumes are positive, proving the last assertion.
\end{proof}

We write \(d=\deg F\). Notice that no realisation of \(\rho_H\) has been assumed. In fact, the finite covering has only been used to calculate \(d\). The same calculation is valid for any discrete cocompact image, with or without torsion.

\subsection{Rong's factorisation}\label{subsec:rong} In the present section we introduce a fundamental building blocks to prove our main result. For this purpose, we first recall the topological description of the vertical pinches and branched coverings as appearing in Rong's Theorem~\ref{thm:rong}. 

Let \(\pi\colon M\to S\) be the projection of our principal circle bundle. For a simple closed curve \(\gamma\subset S\), the restricted circle bundle has total space \(T=\pi^{-1}(\gamma)\), which is an embedded torus formed by the fibres over \(\gamma\). Such a torus is called \textit{vertical}, and it separates \(M\) if and only if \(\gamma\) is a separating curve on \(S\). Let \(\gamma\) be a separating curve and let \(W\) be the closure of one component of \(M\setminus T\). Suppose that \(W\) admits a map \(q\colon W\to V\) to a solid torus such that \(q|_{\partial W}\colon\partial W\to\partial V\) is a homeomorphism. We replace \(W\) by \(V\), using this boundary homeomorphism for the gluing. The resulting degree-one map agrees with \(q\) on \(W\) and with the identity on the retained region and it is called a \textit{vertical pinch}. The new manifold carries a Seifert fibration extending the original one on the retained region. In \S\ref{subsec:pinches}, we show that this induced fibration is again principal. We shall use the term \emph{branched covering} in the topological sense, with the usual codimension-two local model \((z,t)\longmapsto(z^r,t)\), where \(r\ge1\).  In what follows, we shall use the following closed-case form of \cite[Proposition 3.1]{Rong1993}.

\begin{thm}[Rong]\label{thm:rong}
    Let \(M\) and \(N\) be closed, connected, oriented, aspherical Seifert fibred three-manifolds, with fixed Seifert fibrations and infinite fundamental groups. A fibre-preserving map \(F\colon M\to N\) of non-zero degree is homotopic to a composition
    \begin{equation}\label{eq:rong-factorisation}
        M\xrightarrow{\ \Pi\ }P\xrightarrow{\ B\ }N,
    \end{equation}
    where \(\Pi\) is a finite composition of vertical pinches and \(B\) is a fibre-preserving branched covering, branched over fibres. The Seifert fibration of \(P\) is the one induced by the pinches.
\end{thm}

For the map \(B\) appearing in Rong's factorisation, both the branch locus and the branch-value locus are finite unions of fibres. Before moving on, some Remarks are in order.

\begin{rmk}[Comparison between the topological and geometric degrees]\label{rmk:topological-geometric-degree} 
    The non-vanishing assumption in \cite[Proposition~3.1]{Rong1993} concerns the geometric degree, denoted by \(G(F)\), rather than the topological degree \(\deg(F)\). We recall the distinction and explain why the positivity of the latter is sufficient for our application. Since \(M\) and \(N\) are closed, connected and oriented \(3\)-manifolds, the topological degree is characterised by 
    \begin{equation} 
        F_*[M]=\deg(F)\,[N]. 
    \end{equation} 
    In the present setting, Rong's geometric degree is defined as follows. Consider a map \(F'\colon M\to N\) homotopic to \(F\) and an embedded \(3\)-ball \(B\subset N\) such that 
    \begin{equation} 
        (F')^{-1}(B)=B_1\sqcup\cdots\sqcup B_r, 
    \end{equation} 
    where each \(B_j\) is a \(3\)-ball and every restriction 
    \begin{equation} 
        F'|_{B_j}\colon B_j\longrightarrow B 
    \end{equation} 
    is a homeomorphism. The geometric degree \(G(F)\) is the smallest non-negative integer \(r\) obtainable by varying \(F'\) and \(B\); see \cite[Section~1]{Rong1993}. These balls are counted without signs, and the homeomorphisms are allowed to preserve or reverse orientation. For any such choice, let \(\varepsilon_j\in\{+1,-1\}\) record whether \(F'|_{B_j}\) preserves or reverses the orientations induced by \(M\) and \(N\). Homotopy invariance and the local degree formula give 
    \begin{equation} 
        \deg(F)=\deg(F')=\sum_{j=1}^{r}\varepsilon_j. 
    \end{equation} Consequently, \(\lvert\deg(F)\rvert\leq r\). Since this holds for every choice entering the definition of \(G(F)\), taking the minimum yields 
    \begin{equation} 
        \lvert\deg(F)\rvert\leq G(F). 
    \end{equation} 
    In particular, non-zero topological degree implies non-zero geometric degree. In \S\ref{subsec:fundamental-groups-degree}, we established the positivity of \(\deg(F)\), therefore the degree hypothesis required in Rong's proposition holds. 
\end{rmk}

All these hypotheses hold for the mapping \eqref{eq:map-F}: the manifolds are closed, oriented and aspherical with infinite fundamental groups and their fibrations were described in \S\ref{subsec:circle-bundles}, and Lemma \ref{lem:F-surjective-degree} gives \(\deg F=d>0\). As a consequence, the splitting in Equation \eqref{eq:rong-factorisation} holds and there is a homotopy
\begin{equation}\label{eq:F-homotopy}
    F\simeq B\circ\Pi.
\end{equation}
We orient every intermediate manifold using the retained part of its predecessor. Each vertical pinch then has degree \(+1\), so
\begin{equation}\label{eq:degree-B}
    \deg B=d>0.
\end{equation}

\medskip

\subsection{Vertical pinches and surface markings}\label{subsec:pinches} We examine the vertical pinches in Rong's factorisation and their effect on the base surfaces. Each vertical pinch replaces a region, say \(H\times\sph^1\), where \(H\) is a compact subsurface with one boundary component, by a solid torus. We show that the induced fibration on the inserted solid torus is a product fibration \(D^2\times\sph^1\). Consequently, the intermediate manifold remains a principal circle bundle, and the corresponding operation on the base surface replaces \(H\) by a disc, removing its handles. These surface pinches compose to a map, say \(p\colon S\to T\), and we shall verify that \(\Pi_*\) induces precisely \(p_*\) on the fundamental groups of the base surfaces. We then keep track of the marking of \(S\) throughout, so that the resulting factorisation of surface groups retains the prescribed representation \(\rho\). The main result of this subsection is the following.

Before stating the main result of the present section, we recall for the reader's convenience that a pinching map of surfaces is defined as follows. Let \(H\subset S\) be a connected subsurface with one boundary component, and let \(R=\overline{S\setminus H}\). Thus
\begin{equation}
    S=R\cup_{\partial H}H.
\end{equation}
The pinch replaces \(H\) by a disc attached along \(\partial H\), producing
\begin{equation}
    T=R\cup_{\partial H}D^2.
\end{equation}
The resulting map \(p\colon S\to T\) has degree one, agrees with the identity on \(R\), and maps \(H\) to the attached disc while restricting to a homeomorphism on the common boundary. On fundamental groups, \(p_*\) kills the handle generators carried by \(H\), together with their normal closure, and agrees with the inclusion-induced map on \(\pi_1(R)\). We shall also use the term \emph{pinching} for a finite composition of such maps.

\begin{lem}\label{lem:principal-pinch}
    Let \(Q\) be an oriented circle bundle over a closed oriented surface \(S\), and let \(\pi\colon Q\to P\) be a vertical pinch. The induced Seifert fibration on \(P\) has no exceptional fibres. Its base, say \(T\) is obtained from \(S\) by a pinching \(p\colon S\to T\) of surfaces. With compatible orientations of regular fibres,
    \begin{equation}\label{eq:pinch-fibre}
        \pi_*(h_Q)=h_{},
    \end{equation}
    and the homomorphism induced by \(\pi_*\) after quotienting by the fibre subgroups is \(p_*\).
\end{lem}

\begin{rmk}\label{rmk:pinches-exceptional-fibres}
    Recall from Remark~\ref{rmk:source-target-fibrations} that \(M\to S\) is a principal circle bundle, whereas exceptional fibres may occur in \(N\to\mathcal O\). In Rong's factorisation
    \begin{equation}
        M\xrightarrow{\Pi}P\xrightarrow{B}N,
    \end{equation}
    the vertical pinches occur on the source side, starting from \(M\). We shall show that each vertical pinch preserves the absence of exceptional fibres. Consequently, every intermediate manifold arising in the factorisation of \(\Pi\), and in particular \(P\), remains a principal circle bundle over a smooth surface. The exceptional fibres therefore occur only on the target side \(N\), and their interaction with the branched covering \(B\) will be analysed separately.
\end{rmk}

\begin{proof}
Let \(q\colon Q\to S\) be the bundle projection, and let \(\gamma\subset S\) be the separating simple closed curve below the vertical torus along which the pinch is performed. Write
\begin{equation}
    S=S_1\cup_\gamma S_2,
\end{equation}
and set
\begin{equation}
    W_i=q^{-1}(S_i),\qquad i=1,2.
\end{equation}
We assume that \(W_2\) is the region replaced by the pinch, while \(W_1\) is retained. Since \(S_2\) is a compact oriented surface with one boundary component and an oriented circle bundle over a surface with boundary is trivial, we may identify
\begin{equation}\label{eq:pinch-product}
    W_2=S_2\times\sph^1.
\end{equation}

The vertical pinch replaces \(W_2\) by a solid torus \(V\). Write \(f=\pi|_{W_2}\colon W_2\longrightarrow V\), and identify the common boundary \( \partial W_1=\partial W_2\). Let \(\varphi=f|_{\partial W_2}\colon\partial W_2\longrightarrow\partial V\) be the boundary homeomorphism. The manifold \(P\) is obtained by gluing \(V\) to \(W_1\) using \(\varphi\).

In \(H_1(\partial W_2;\Z)\), let \(c\) be the homology class represented by the boundary of a section, say \(S_2\times\{t\}\), let \(h\) be the positively oriented fibre class, and notice that \(\{c,h\}\) form a basis. Since \(S_2\) has only one boundary component, its oriented boundary is null-homologous in \(S_2\). Hence the product description \eqref{eq:pinch-product} gives
\begin{equation}\label{eq:boundary-kernel}
    \ker\bigl(H_1(\partial W_2;\Z)\longrightarrow H_1(W_2;\Z)\bigr)
    =\Z c.
\end{equation}

Let \(i_W\colon\partial W_2\hookrightarrow W_2\) and \(i_V\colon\partial V\hookrightarrow V\) denote the inclusion maps. By construction, the class \(c\) is represented by the boundary of the section \(S_2\times\{t\}\subset W_2\), and hence \((i_W)_*(c)=0\) in \(H_1(W_2;\Z)\). Since \(f|_{\partial W_2}=\varphi\), the diagram
\begin{equation}
    \begin{tikzcd}[column sep=large, row sep=large]
        \partial W_2
        \arrow[r,"\varphi"]
        \arrow[d,"i_W"']
        &
        \partial V
        \arrow[d,"i_V"]
        \\
        W_2
        \arrow[r,"f"']
        &
        V
    \end{tikzcd}
\end{equation}
commutes. Therefore \((i_V)_*\bigl(\varphi_*(c)\bigr)=f_*\bigl((i_W)_*(c)\bigr)=0\), and hence 
\begin{equation}
    \varphi_*(c)\in
    \ker\bigl(H_1(\partial V;\Z)\longrightarrow H_1(V;\Z)\bigr).
\end{equation}

For a solid torus \(V\), this kernel is generated by the meridian class \(\mu\). Fix a longitude class \(\lambda\) such that \(\{\mu,\lambda\}\) is an oriented basis of \(H_1(\partial V;\Z)\). Under the inclusion
\begin{equation}
    H_1(\partial V;\Z)\longrightarrow H_1(V;\Z),
\end{equation}
the meridian \(\mu\) becomes trivial, whereas \(\lambda\) maps to a generator of \(H_1(V;\Z)\cong\Z\). Hence
\begin{equation}
    \ker\bigl(H_1(\partial V;\Z)\longrightarrow H_1(V;\Z)\bigr)=\Z\mu.
\end{equation}

The class \(c\) is primitive in \(H_1(\partial W_2;\Z)\), since it belongs to the basis \(\{c,h\}\). Since \(\varphi\) is a homeomorphism, the induced map \(\varphi_*\colon H_1(\partial W_2;\Z)\longrightarrow H_1(\partial V;\Z)\) is an isomorphism and therefore preserves primitive classes. Thus \(\varphi_*(c)\) is a primitive element of \(\Z\mu\), which forces
\begin{equation}
    \varphi_*(c)=\pm\mu.
\end{equation}

Write \( \varphi_*(h)=a\mu+b\lambda\), where \(a,b\in\Z\). Since \(\{c,h\}\) is a basis of \(H_1(\partial W_2;\Z)\), the matrix of \(\varphi_*\) with respect to the bases \(\{c,h\}\) and \(\{\mu,\lambda\}\) has the form
\begin{equation}
    \begin{pmatrix}
        \pm1 & a\\
        0 & b
    \end{pmatrix}.
\end{equation}
Since \(\varphi_*\) is an isomorphism, this matrix has determinant \(\pm1\). Therefore \(b=\varepsilon\in\{1,-1\}\), and hence
\begin{equation}\label{eq:fibre-longitude}
    \varphi_*(h)=a\mu+\varepsilon\lambda,
    \qquad
    a\in\Z,\quad \varepsilon\in\{1,-1\}.
\end{equation}
Since \(\{\mu,\lambda\}\) is an oriented basis, \(\mu\cdot\lambda=1\), and therefore
\begin{equation}
    \bigl|\mu\cdot\varphi_*(h)\bigr|
    =
    |\varepsilon|
    =
    1.
\end{equation}

In a fibred solid torus, the absolute intersection number between the meridian and the regular fibre slope on the boundary is the multiplicity of the core fibre. Hence the core fibre of \(V\) has multiplicity one and is therefore regular. Thus the filling introduces no exceptional fibre. Since the fibration on the retained region \(W_1\) is unchanged and \(Q\to S\) is a principal circle bundle, all fibres of \(P\) are regular. Consequently, the induced Seifert fibration on \(P\) is again a principal circle bundle. This is also the content of \cite[Lemma~3.2.8]{Konstantinov2019}. The base of the induced fibration on \(V\) is a disc. Therefore the base of \(P\) is
\begin{equation}
    T=S_1\cup_\gamma D^2,
\end{equation}
where the boundary identification is induced by \(\varphi\). The orientation of \(T\) extends that of \(S_1\), and the fibre orientation extends the one on \(W_1\). Let
\begin{equation}
    p\colon S\longrightarrow T
\end{equation}
be the corresponding ordinary surface pinch: it agrees with the retained identification on \(S_1\) and replaces \(S_2\) by the disc \(D^2\). We now identify the homomorphism induced on the fundamental groups of the base surfaces. Choose a base point, say \(s_0\in\gamma\), choose a point \(z_0\in q^{-1}(s_0)\), and use their images as base points in the target spaces. Since \(\pi\) agrees with the identity on the retained region \(W_1\), a positively oriented fibre through \(z_0\) is mapped to a positively oriented fibre of \(P\). Hence
\begin{equation}\label{eq:pinch-fibre2}
    \pi_*(h_Q)=h_P.
\end{equation}

The bundle projections induce natural identifications
\begin{equation}
    \frac{\pi_1(Q)}{\langle h_Q\rangle}\cong\pi_1(S),
    \qquad
    \frac{\pi_1(P)}{\langle h_P\rangle}\cong\pi_1(T).
\end{equation}
Equation~\eqref{eq:pinch-fibre2} therefore implies that \(\pi_*\) descends to a homomorphism \( \overline{\pi}_*\colon\pi_1(S)\longrightarrow\pi_1(T)\). We claim that \(\overline{\pi}_*=p_*\). Choose a section of \(W_2\to S_2\) passing through \(z_0\), and use it to lift based loops in \(S_2\) to \(W_2\). Their images under \(f\) lie in the solid torus \(V\). By \eqref{eq:fibre-longitude}, the image of a boundary fibre represents a generator of
\begin{equation}
    \pi_1(V)\cong\Z,
\end{equation}
since its longitude coefficient is \(\pm1\). Hence the images in \(\pi_1(P)\) of the lifted loops from \(S_2\) lie in the fibre subgroup \(\langle h_P\rangle\). They therefore become trivial after passing to \(\pi_1(T)\). Thus \(\overline{\pi}_*\) kills the subgroup of \(\pi_1(S)\) carried by \(S_2\). On the other hand, since \(\pi\) agrees with the identity on \(W_1\), the restriction of \(\overline{\pi}_*\) to loops contained in \(S_1\) agrees with the homomorphism induced by the inclusion
\begin{equation}
    S_1\hookrightarrow T.
\end{equation}
These are precisely the two defining properties of \(p_*\): the surface pinch \(p\) kills the handle generators carried by \(S_2\) and agrees with the retained identification on \(S_1\). Van Kampen's theorem therefore gives
\begin{equation}
    \overline{\pi}_*=p_*,
\end{equation}
as required.
\end{proof}

Iterating Lemma \ref{lem:principal-pinch} shows that \(P\) is an oriented circle bundle \(q_P\colon P\longrightarrow T\) over a smooth closed oriented surface \(T\). There is a composition of surface pinches \(p\colon S\to T\) for which
\begin{equation}\label{eq:Pi-quotient}
    \Pi_*(h_M)=h_P,\qquad \overline{\Pi}_*=p_*;
\end{equation}
where the bar denotes the map on the quotients by the fibre subgroups. The assertion is about induced homomorphisms; it does not require the pinch maps themselves to commute pointwise with the bundle projections. By combining identities \eqref{eq:F-homotopy} and \eqref{eq:Pi-quotient} with the identity \(F_*(h_M)=h_N\), we obtain
\begin{equation}
\label{eq:B-fibre-degree}
B_*(h_P)=h_N.
\end{equation}
Base-point changes introduce conjugations only, and do not affect this identity because the fibre subgroups are central. In particular, \(B\) has signed regular-fibre degree \(+1\). By quotienting \(B_*\) by the fibre subgroups gives a homomorphism
\begin{equation}
\label{eq:beta-representation}
\beta\colon\pi_1(T)\longrightarrow G\hookrightarrow\pslr2.
\end{equation}
The same quotient applied to \eqref{eq:F-homotopy} yields the marked identity
\begin{equation}
\label{eq:rho-factorisation}
\rho=\beta\circ p_*,
\end{equation}
up to a single overall conjugation. Since \(F_*=B_*\circ\Pi_*\) is surjective, \(B_*\) is surjective as well. Passing to the quotients by the fibre subgroups, we conclude that \(\beta\colon\pi_1(T)\to G\) is onto. Since \(G\) is a non-abelian cocompact hyperbolic orbifold group, \(\pi_1(T)\) is non-abelian. As a consequence, \(T\) has genus at least two and is therefore aspherical.

\smallskip

\subsection{Local degrees and cone angles}\label{subsec:local-degrees}

The fibre-preserving map \(B\) induces a map of the underlying surfaces \(b\colon T\longrightarrow|\,\mathcal O\,|\). We shall verify that pulling back the orbifold metric by \(b\) gives integral cone angles. For this purpose, we first provide a description of local degrees and cone angles. Recall from \S\ref{sssec:fibre-preserving-maps} that, when a regular source fibre is mapped to an exceptional target fibre of multiplicity \(m\), the local degree of the induced map on the underlying base surfaces contains two contributions; namely the transverse ramification index \(r\) of \(B\) and the multiplicity \(m\) of the target fibre. Thus the local degree of the base map is \(mr\), whereas the genuine transverse ramification of the mapping \(B\) is only \(r\).

\begin{lem}\label{lem:local-degree}
    Let \(x\in T\), and suppose that \(b(x)\) has orbifold order \(m\geq1\), with \(m=1\) at a regular point. If the transverse ramification index of \(B\) along the fibre over \(x\) is equal to \(r\), then the local degree of \(b\) at \(x\) is \(mr\). Moreover, the pull-back cone angle at \(x\) is \(2\pi r\).
\end{lem}

\begin{proof}
Since \(P\to T\) is a principal circle bundle, the fibre over \(x\) is regular, while the fibre over \(b(x)\) has multiplicity \(m\). Moreover, by \eqref{eq:B-fibre-degree}, \(B\) has degree \(+1\) on regular fibres. The local description of \S\ref{sssec:fibre-preserving-maps} therefore applies and gives
\begin{equation}
    \deg_x(\,b\,)=mr,
\end{equation}
while the pull-back of the orbifold angle \(2\pi/m\) is
\begin{equation}\label{eq:pullback-angle}
    mr\cdot\frac{2\pi}{m}=2\pi r,
\end{equation}
as desired.
\end{proof}

Away from the finitely many branch values and cone points, \(b\) is an orientation-preserving covering and its total degree is \(d\). In fact, recall that \(d=\deg F\) and that Rong's factorisation gives \(\deg B=d\). Since \(B\) has degree \(+1\) on regular fibres, the induced map on the base surfaces has degree \( \deg b=\deg B=d\) as desired. Now, let \(y\in\mathcal O\) be a cone point of order \(m\). By Lemma~\ref{lem:local-degree}, every \(x\in b^{-1}(y)\) has local degree \(mr_x\), where \(r_x\) is the transverse ramification index of \(B\) along the fibre over \(x\). Hence
\begin{equation}\label{eq:cone-degree-divisibility}
    d
    =\sum_{x\in b^{-1}(y)}\deg_x b
    =m\sum_{x\in b^{-1}(y)}r_x.
\end{equation}
In particular, every cone order of \(\mathcal O\) divides \(d\). Since \(G\) contains torsion, \(\mathcal O\) has at least one cone point of order \(m\geq2\), and therefore \(d\geq2\). Moreover, the local model of \S\ref{sssec:fibre-preserving-maps} shows that the pull-back of the hyperbolic orbifold metric defines a branched hyperbolic structure on \(T\), with cone angle \(2\pi r_x\) at \(x\). It remains to identify its holonomy.

\smallskip

\subsection{The developing map and its holonomy}\label{subsec:holonomy-identification} As already alluded above, the fibre-preserving map \(B\) induces a map of the underlying surfaces \(b\colon T\to |\mathcal O|\). By the previous discussion, pulling back the hyperbolic orbifold metric by \(b\) produces a branched hyperbolic metric on \(T\) with integral cone angles. It remains to show that this metric is induced by a developing map whose holonomy is precisely the representation \(\beta\) defined above.

\begin{lem}\label{lem:base-holonomy}
    The pull-back of the orbifold metric under \(b\) is a branched hyperbolic structure \(\sigma_T\) on \(T\), and its holonomy is the representation \(\beta\) in \eqref{eq:beta-representation}.
\end{lem}

\begin{proof}
Let \(\pi_T\colon\widetilde T\longrightarrow T\) be the universal covering projection, and let \( \widehat P\longrightarrow \widetilde T\) be the pull-back of the circle bundle \(q_P\colon P\longrightarrow T\). Recall that this defined as \(\widehat P=\bigl\{(x,z)\in\widetilde T\times P\mid\pi_T(x)=q_P(z)\bigr\}\). Denote by \(\widehat q_P\colon\widehat P\longrightarrow\widetilde T\) and \( r_P\colon\widehat P\longrightarrow P \) the two natural projections. Thus the diagram
\begin{equation}\label{eq:pullback-P-diagram}
    \begin{tikzcd}[column sep=large, row sep=large]
        \widehat P
        \arrow[r,"r_P"]
        \arrow[d,"\widehat q_P"']
        &
        P
        \arrow[d,"q_P"]
        \\
        \widetilde T
        \arrow[r,"\pi_T"']
        &
        T
    \end{tikzcd}
\end{equation}
is commutative. Since \(T\) has genus at least two, its universal cover \(\widetilde T\) is contractible. Consequently, the pull-back circle bundle is trivial, that is, \(\widehat P\cong\widetilde T\times\sph^1\). In particular, \(\pi_1(\widehat P)\cong\Z\). We next observe that \(r_P\colon\widehat P\to P\) is the covering associated with the kernel of \((q_P)_*\colon\pi_1(P)\longrightarrow\pi_1(T)\). By the homotopy exact sequence of the circle bundle, this kernel is precisely the subgroup generated by the positively oriented fibre \(h_P\). As a consequence, \((r_P)_*\bigl(\pi_1(\widehat P)\bigr)=\langle h_P\rangle\). The regular covering
\begin{equation}\label{eq:tangent-G-cover}
    q\colon T^1\HH^2\longrightarrow N
\end{equation}
corresponds to the subgroup \(\langle h_N\rangle\) of \(\pi_1(N)\). Since \(B_*(h_P)=h_N\) by \eqref{eq:B-fibre-degree}, it follows that
\begin{equation}
    (B\circ r_P)_*\bigl(\pi_1(\widehat P)\bigr)
    \subset
    \langle h_N\rangle
    =
    q_*\bigl(\pi_1(T^1\HH^2)\bigr).
\end{equation}
The covering-space lifting criterion therefore gives a lift \(\widehat B\colon\widehat P\longrightarrow T^1\HH^2\), such that
\begin{equation}\label{eq:lifted-B}
    q\circ\widehat B=B\circ r_P.
\end{equation}
The deck groups of \(r_P\) and \(q\) are naturally identified with
\begin{equation}
    \pi_1(P)/\langle h_P\rangle\cong\pi_1(T)
    \qquad\textnormal{and}\qquad
    \pi_1(N)/\langle h_N\rangle\cong G,
\end{equation}
respectively. Since \(B_*(h_P)=h_N\), the homomorphism \(B_*\) induces a homomorphism between these quotient groups, which is precisely the representation \(\beta\colon\pi_1(T)\longrightarrow G\). With compatible choices of lifts and base points, we obtain
\begin{equation}\label{eq:lifted-B-equivariance}
    \widehat B(\gamma z)=\beta(\gamma)\widehat B(z) \qquad \textnormal{for }\gamma\in\pi_1(T).
\end{equation}


Let \(\operatorname{pr}\colon T^1\HH^2\to\HH^2\) be the projection, let \(q_N\colon N\to\mathcal O\) be the Seifert projection and let \(\pi_{\mathcal O}\colon\HH^2\to\mathcal O\) be the quotient map. The projections satisfy
\begin{equation}
    q_N\circ q
    =
    \pi_{\mathcal O}\circ\operatorname{pr}.
\end{equation}
We claim that \(\operatorname{pr}\circ\widehat B\) is constant on every circle fibre of \(\widehat P\to\widetilde T\). Let \(C\subset\widehat P\) be such a fibre. The map \(r_P\) sends \(C\) onto a fibre of \(P\to T\), and, since \(B\) is fibre preserving, \(B(r_P(C))\) is contained in a single Seifert fibre of \(N\to\mathcal O\). By using \(q\circ\widehat B=B\circ r_P\), we obtain
\begin{equation}
    \pi_{\mathcal O}\bigl(\operatorname{pr}(\widehat B(C))\bigr)
    =
    q_N\bigl(B(r_P(C))\bigr),
\end{equation}
and the right-hand side is a single point of \(\mathcal O\). Therefore
\(\operatorname{pr}(\widehat B(C))\) is contained in a single \(G\)-orbit in \(\HH^2\). Since \(G\) is discrete, this orbit is a discrete subset of \(\HH^2\). The circle \(C\) is connected, and hence the continuous image
\(\operatorname{pr}(\widehat B(C))\) consists of a single point. Thus
\(\operatorname{pr}\circ\widehat B\) is constant on every circle fibre of \(\widehat P\). 
It follows that \(\operatorname{pr}\circ\widehat B\) descends to a map
\begin{equation}
\label{eq:dev-T}
    \textnormal{dev}_T\colon \widetilde T\longrightarrow\HH^2,
\end{equation}
and the identity \eqref{eq:lifted-B-equivariance} implies
\begin{equation}\label{eq:dev-T-equivariance}
    \textnormal{dev}_T(\gamma x)=\beta(\gamma)\,\textnormal{dev}_T(x).
\end{equation}

Locally, this map is the transverse lift described in Lemma \ref{lem:local-degree}, possibly followed by an element of \(G\). Hence it is an orientation-preserving topological branched map of local degree \(r\). We explain how it determines the required complex structure, without assuming one on \(T\). At a point of \(\widetilde T\), choose a conformal coordinate \(z\) about its image in \(\HH^2\), centred at that image. Classification of coverings of punctured discs gives an oriented source coordinate \(w\), extending continuously over the centre, such that
\begin{equation}\label{eq:complex-branch-chart}
    z\circ\textnormal{dev}_T=w^r.
\end{equation}
Away from the ramification points these are pull-backs of conformal charts, so their transition maps are holomorphic. Across a ramification point, each transition map and its inverse are continuous and holomorphic on a punctured neighbourhood; removable singularities make them holomorphic at the centre as well. These charts therefore define a complex structure on \(\widetilde T\). Equivariance and the holomorphic action of \(\beta(\pi_1(T))\) imply that its deck transformations are holomorphic, first off the ramification points and then everywhere by the same extension argument. The complex structure consequently descends to \(T\), and \(\textnormal{dev}_T\) is holomorphic. Finally, the metric \(\textnormal{dev}_T^*h_{\HH}\) descends to \(T\) by \eqref{eq:dev-T-equivariance}. It agrees away from the finite exceptional set with the pull-back of the orbifold metric by \(b\), and it has cone angles \(2\pi r\) there. This proves the first assertion. Equation \eqref{eq:dev-T-equivariance} identifies its holonomy with the prescribed homomorphism \(\beta\), proving the second.
\end{proof}

\subsection{Existence of genuine ramification}
\label{subsec:genuine-ramification}

The next point is needed to reverse the surface pinches. We must know that \(\sigma_T\) has at least one branch point of positive order.

\begin{lem}
\label{lem:genuine-ramification}
The branched covering \(B\) has non-empty ramification locus. Consequently, \(\sigma_T\) has a point of cone angle at least \(4\pi\).
\end{lem}

\begin{proof}
Suppose that \(B\) has no ramification. Then its local degree is one everywhere, so \(B\) is a local homeomorphism. Since \(P\) is compact, \(B\) is proper and hence is an ordinary covering of \(N\). As the source is connected, the number of sheets of this covering is equal to \(\bigl[\pi_1(N):B_*(\pi_1(P))\bigr]\). Since \(B_*\) is surjective, this index is one. Therefore \(B\) is a one-sheeted covering and hence a homeomorphism. By \eqref{eq:B-fibre-degree}, this homeomorphism sends the regular-fibre subgroup of \(\pi_1(P)\) onto that of \(\pi_1(N)\). Passing to the quotient gives
\begin{equation}
\label{eq:torsion-contradiction}
\pi_1(T)
\cong\frac{\pi_1(P)}{\langle\, h_P\,\rangle}
\cong\frac{\pi_1(N)}{\langle\, h_N\,\rangle}
\cong G.
\end{equation}
The left-hand group is torsion-free, whereas \(G\) contains torsion. This is a contradiction. Thus some fibre of \(P\) has transverse ramification index \(r\geq2\). By Lemma \ref{lem:local-degree}, the corresponding point of \(\sigma_T\) has angle \(2\pi r\geq4\pi\).
\end{proof}

The transverse index is the relevant one here. A regular fibre can cover an exceptional fibre several times even when the map of three-manifolds is unbranched. That multiplicity alone would not produce a branch point on \(T\).

\subsection{Reinserting the pinched handles}\label{subsec:handle-insertion}

We now give the local construction that reverses the surface pinches while preserving the marked holonomy. It is the handle-insertion operation used in \cite[Lemma 3.12]{Faraco2021}. We give a construction by monodromy to explain why the argument applies to arbitrary holonomy and to any number of removed handles.

\begin{lem}\label{lem:handle-insertion}
    Let \(\sigma_T\) be a branched hyperbolic structure on a closed oriented surface \(T\), with holonomy \(\beta\), and suppose that \(\sigma_T\) has a branch point. If \(p\colon S\to T\) is a composition of ordinary surface pinches, then \(\beta\circ p_*\) is geometrisable on \(S\).
\end{lem}

\begin{proof}
    We begin by realising the replacement for a single local branched chart. Choose a closed disc \(U\) about a branch point of local degree \(r\geq2\), contained in a slightly larger branched chart. Its developing map restricts to an orientation-preserving branched covering, say \(d_U\colon U\longrightarrow V\) of degree \(r\), where \(V\subset\HH^2\) is a closed disc and the only branch value is in its interior. Near the boundary, \(d_U\) is a connected covering of annuli of degree \(r\). Fix an integer \(k\geq1\). The idea is to replace \(U\) by a compact oriented surface \(W_{k,1}\) of genus \(k\) with one boundary component, equipped with a degree-\(r\) branched covering, say
    \begin{equation}\label{eq:new-disc-map}
        q_k\colon W_{k,1}\longrightarrow V
    \end{equation}
    whose boundary covering agrees with that of \(d_U\). Choose \(r-1+2k\) distinct interior points of \(V\). On the \(r\) sheets, prescribe the monodromy about these points, in an ordered system of peripheral loops, by
    \begin{equation}\label{eq:disc-monodromy}
        (1\,2),(2\,3),\ldots,(r-1\,r),
            \underbrace{(1\,2),(1\,2),\ldots,(1\,2),(1\,2)}_{2k\ \textnormal{terms}}.
    \end{equation}
Let \(V^\circ\) be the disc obtained from \(V\) by removing the chosen branch values. Its fundamental group is freely generated by peripheral loops around these punctures. Assigning to these generators the transpositions listed in \eqref{eq:disc-monodromy} therefore defines a monodromy representation \(\pi_1(V^\circ)\longrightarrow\mathfrak S_r\), and hence an \(r\)-sheeted covering of \(V^\circ\). Since the adjacent transpositions \((1\,2),(2\,3),\ldots,(r-1\,r)\) generate a transitive subgroup of \(\mathfrak S_r\), the covering is connected. Around each puncture the local monodromy is a transposition, so the covering extends across that puncture by inserting one simple branch point together with \(r-2\) regular preimages. Filling all punctures in this way yields a compact connected oriented branched covering of \(V\).

Notice that the product of these transpositions is an \(r\)-cycle, up to the convention used for composing permutations. Hence the total boundary monodromy is again an \(r\)-cycle, or its inverse according to the orientation convention for the boundary loop. In either case the boundary of the covering surface is connected and has degree \(r\) over \(\partial V\). Riemann--Hurwitz, or an Euler characteristic calculation after deleting the branch values, gives
    \begin{equation}\label{eq:replacement-euler}
        \chi(W_{k,1})
        =r\chi(V)-(r-1+2k)
        =1-2k.
    \end{equation}
    The covering surface thus has genus \(k\) and one boundary component, as required. Choose an annular collar \(C\) of \(\partial V\) containing no branch value. Both \(d_U^{-1}(C)\) and \(q_k^{-1}(C)\) are connected degree-\(r\) coverings of \(C\). Classification of connected coverings of an annulus gives an orientation-preserving isomorphism
    \begin{equation}\label{eq:collar-identification}
        \psi\colon q_k^{-1}(C)\longrightarrow d_U^{-1}(C), \qquad d_U\circ\psi=q_k.
    \end{equation}
    Pull back the hyperbolic metric of \(V\) under \(q_k\). Equation \eqref{eq:collar-identification} identifies this metric and its developing charts with the old ones on a whole collar, so replacing \(U\) by \(W_{k,1}\) produces a branched hyperbolic structure. The new patch has a global developing map \(q_k\) with values in \(V\subset\HH^2\). By construction, the holonomy of every loop in \(W_{k,1}\) is trivial. 
    
    \begin{rmk}\label{rmk:prescribed-pinch}
        The local construction above inserts the genus-\(k\) subsurface at the chosen branch point. We next need to compare this geometric realisation with an arbitrary prescribed pinching map \(p\colon S\to T\). Suppose that \(p\) collapses a genus-\(k\) subsurface \(A\subset S\), with one boundary component, onto a disc \(D\subset T\), whereas the local handle-insertion construction is performed in a disc \(U\subset T\) containing a branch point of \(\sigma_T\). In general, the two discs need not coincide. Choose the base point outside \(D\cup U\). There exists an orientation-preserving ambient isotopy of \(T\), fixed at the base point, whose final map \(a\colon T\to T\) sends \(D\) onto \(U\). Replacing \(p\) by
        \begin{equation}
            p'=a\circ p
        \end{equation}
    therefore moves the disc associated with the pinch to the disc where the local construction is performed. Since \(a\) is isotopic to the identity relative to the base point, then the mapping \(p'\) is based homotopic to \(p\), and hence
        \begin{equation}\label{eq:marked-pinch-isotopy}
            p'_* = p_*.
        \end{equation}
    Thus, without changing the homomorphism induced by the prescribed pinching map, we may assume that the disc produced by the pinch is precisely the disc \(U\) containing the chosen branch point.
\end{rmk}

    We now complete the construction for the prescribed marking. Suppose first that \(p\colon S\to T\) is a single pinch collapsing a genus-\(k\) subsurface \(H\subset S\), with one boundary component. By Remark~\ref{rmk:prescribed-pinch}, after replacing \(p\) within its based homotopy class, we may assume that \(p\) maps the subsurface \(H\) onto the disc \(U\) used in the local construction above. We take \(p\) in standard pinch form: it identifies the retained subsurface with \(T\setminus\operatorname{int}(U)\) and restricts to a homeomorphism on a collar of the common boundary. We now identify \(H\) with the replacement surface \(W_{k,1}\). The collar identification in \eqref{eq:collar-identification} provides an orientation-preserving identification between neighbourhoods of \(\partial W_{k,1}\) and \(\partial U\), compatible with the developing maps. Since \(H\) and \(W_{k,1}\) are compact oriented surfaces of genus \(k\) with one boundary component, this boundary identification extends to an orientation-preserving homeomorphism
    \begin{equation}
        H\longrightarrow W_{k,1}.
    \end{equation}
    By using this identification, we keep the original branched hyperbolic structure on the retained subsurface and equip \(H\cong W_{k,1}\) with the structure induced by \(q_k\). Since the developing maps agree on an entire collar of the common boundary, the two structures glue to a branched hyperbolic structure \(\sigma_S\) on \(S\). It remains to identify its holonomy. On the retained subsurface, the new structure agrees, through the identification induced by \(p\), with the restriction of \(\sigma_T\) to \(T\setminus\operatorname{int}(U)\). Its holonomy is therefore the restriction of \(\beta\circ p_*\). On \(H\cong W_{k,1}\), the map
    \begin{equation}
        q_k\colon W_{k,1}\longrightarrow V\subset\HH^2
    \end{equation}
    is a globally defined developing map, and hence every loop contained in \(H\) has trivial holonomy. This agrees with \(\beta\circ p_*\), since the pinch kills the handle generators carried by \(H\). Finally, since the subgroups carried by the retained subsurface and by \(H\), together with the common boundary relation, generate \(\pi_1(S)\), van Kampen's theorem gives
    \begin{equation}\label{eq:inserted-holonomy}
        \hol(\sigma_S)=[\beta\circ p_*].
    \end{equation}

    Finally, suppose that \(p\) is a composition of pinches. Write
    \begin{equation}
        S=S_0\xrightarrow{\,p_1\,}S_1\xrightarrow{\,p_2\,}\cdots\xrightarrow{\,p_\ell\,}S_\ell=T,
    \end{equation}
    where \(p_i\) collapses a genus-\(k_i\) subsurface with one boundary component. We apply the construction in reverse order. At the \(i\)-th step, replacing a disc with surface \(W_{k_i,1}\) decreases the Euler characteristic by \(2k_i\), and hence increases the genus by \(k_i\). The resulting surface is therefore orientation-preserving homeomorphic to \(S_{i-1}\). After all the insertions, the final surface is orientation-preserving homeomorphic to \(S_0=S\). The replacement described in \eqref{eq:disc-monodromy} contains simple branch points, so after each insertion a branch point remains available for the next one. Repeating the marking and holonomy argument at every step yields a branched hyperbolic structure on \(S\) with holonomy \([\beta\circ p_*]\), completing the proof.
\end{proof}

The numerical effect agrees with the Euler-number formula. Replacing \(U\) by \(W_{k,1}\) decreases the Euler characteristic by \(2k\). The old branching order in \(U\) was \(r-1\); the new patch has \(r-1+2k\) simple branch points, so the total branching order increases by \(2k\). Thus \(\chi+\kappa\), and hence the holonomy Euler number, remains unchanged. The area is unchanged as well, since both patches have area \(r\,\textnormal{Area}(V)\). The assumption \(r\geq2\) is exactly what permits the transpositions in the construction.

\subsection{Completion of the proof}\label{subsec:completion} We are finally in the right position to complete the proof of Theorem \ref{thm:discrete-geometrisation}. We begin with the following

\begin{proof}[Proof of Proposition \ref{prop:torsion-geometrisation}]
    The argument developed in \S\S\ref{subsec:circle-bundles}--\ref{subsec:pinches} provide a pinching map of surfaces, that is, \(p\colon S\to T\) and a representation \(\beta\colon\pi_1(T)\to G\) such that \(\rho=\beta\circ p_*\). By Lemma \ref{lem:base-holonomy}, \(\beta\) is the holonomy of the branched hyperbolic structure \(\sigma_T\). Lemma \ref{lem:genuine-ramification} guarantees the existence of a branch point and, finally, Lemma \ref{lem:handle-insertion} provides a branched hyperbolic structure \(\sigma\) on \(S\) with the desired holonomy \([\,\rho\,]\). If the factorisation has no pinches, the structure \(\sigma_T\) already gives the desired realisation. This proves the proposition.
\end{proof}

\begin{proof}[Proof of Theorem \ref{thm:discrete-geometrisation}]
By Lemmata \ref{lem:volume-euler} and \ref{lem:noncompact-euler-zero}, positive volume implies that \(\mathcal O=G\backslash\HH^2\) is compact. We first prove necessity of the stated condition. Suppose that \(\rho\) is the holonomy of a branched hyperbolic structure on \(S\). Its developing map descends to a map \(f\colon S\longrightarrow|\,\mathcal O\,|\). At a point whose developing map has local degree \(r\), and whose image in \(\mathcal O\) has orbifold order \(m\), the local degree of \(f\) is \(mr\). Thus \(f\) is an orientation-preserving branched covering of the underlying surfaces. Its image is open and compact, hence is all of \(|\,\mathcal O\,|\). Integrating the pull-back orbifold area form away from its finitely many exceptional values gives
\begin{equation}\label{eq:necessary-degree}
    \textnormal{vol}(\,\rho\,)=\deg(\,f\,)\,\textnormal{vol}(\,\mathcal O\,),   \qquad  \deg(\,f\,)\in\Z_{>0}.
\end{equation}

If the two volumes are equal, then \(\deg(f)=1\). Every local degree \(mr\) is therefore one. Consequently, \(\mathcal O\) has no cone points, the structure has no branch points, and \(f\) is a homeomorphism. The holonomy is the composition of the isomorphism \(f_*\) with the Fuchsian inclusion \(\pi_1(\mathcal O)\hookrightarrow\pslr2\), and is faithful. This proves that a geometrisable representation satisfies the alternative in the theorem. 

For sufficiency, let us first suppose that \(G\) contains torsion. Proposition \ref{prop:torsion-geometrisation} supplies the desired structure and Equation \eqref{eq:cone-degree-divisibility} also shows directly that the strict volume inequality holds. Suppose that \(G\) is torsion-free. The quotient \(\Sigma=G\backslash\HH^2\) is a closed surface equipped with a complete hyperbolic structure and the representation \(\rho\colon\pi_1(S)\longrightarrow \pi_1(\Sigma)\cong G<\pslr2\) is induced by a map \(f\colon S\longrightarrow\Sigma\). Such a map exists and is unique up to based homotopy because \(\Sigma\) is aspherical. Naturality gives
\begin{equation}\label{eq:torsion-free-degree}
    \deg(f) =\frac{e(\,\rho\,)}{\chi(\,\Sigma\,)} =\frac{\textnormal{vol}(\,\rho\,)}{\textnormal{vol}(\,\Sigma\,)} >0.
\end{equation}
In particular, \(\deg(f)\) is a positive integer. If \(\rho\) is faithful, then \(f_*\) is an isomorphism. As both surfaces are aspherical, \(f\) is a homotopy equivalence, and its positive degree is one. A degree-one map of closed oriented surfaces is homotopic to a pinch followed by an orientation-preserving homeomorphism, by Edmonds' factorisation theorem \cite{Edmonds1979}. The injectivity of \(f_*\) rules out every non-trivial pinch. Thus \(f\) is homotopic to an orientation-preserving homeomorphism, and pulling back the hyperbolic metric of \(\Sigma\) gives the required unbranched structure. If instead \eqref{eq:main-volume-condition} holds, then \(\deg(f)>1\). In particular, \(f\) is not homotopic to a pinch map, since a pinch has degree one. The representation is non-elementary and purely hyperbolic, and its image is cocompact. It is also non-Fuchsian: faithfulness would make \(f\) a homotopy equivalence of degree one, as above. The geometrisation criterion in \cite[Theorem 3.14]{Faraco2021}, applied in this positive-degree setting, therefore gives a branched hyperbolic structure with holonomy \([\rho]\). This completes the proof.
\end{proof}

For every discrete representation of positive volume, the ratio \(\textnormal{vol}(\rho)/\textnormal{vol}(\mathcal O)\) is a positive integer, by \eqref{eq:degree-formula} and \eqref{eq:torsion-free-degree}. Hence either the two volumes are equal or \(\textnormal{vol}(\rho)>\textnormal{vol}(\mathcal O)\). It follows that the only discrete representations of positive volume which are not geometrisable are the non-faithful ones for which \(\textnormal{vol}(\rho)=\textnormal{vol}(\mathcal O)\). The theorem proves the discrete case of the strong form of Goldman's conjecture stated in \cite[Conjecture 1]{DeroinTholozan2026}. 

\medskip
\section{Symmetric branched hyperbolic structures}\label{sec:sym structures}

The main theorem has a natural consequence for the realisation of prescribed symmetries on finite covers. Indeed, passing to a normal torsion-free finite-index subgroup produces a regular covering whose deck group is prescribed algebraically, and the geometrisation result allows us to realise the restricted holonomy by a branched hyperbolic structure invariant under this deck action. 

\begin{thm}\label{cor:invariant-cover}
    Let \(\rho\colon\pi_1(S)\to\pslr2\) have discrete cocompact image \(G\). 
    Let \(H\trianglelefteq G\) be torsion-free and of finite index, and let \(p_H\colon S_H\longrightarrow S\) be the covering corresponding to \(\rho^{-1}(\,H\,)\), with the lifted orientation. Suppose \(\rho_H=\rho|_{\pi_1(\,S_H\,)}\) is geometrisable, then it admits a realisation invariant under the prescribed deck action of \(G/H\) on \(S_H\). This realisation is the pull-back of a branched hyperbolic structure on \(S\) with holonomy \(\rho\). In particular, \(G/H\) acts by orientation-preserving isometries of the resulting branched hyperbolic structure on \(S_H\).
\end{thm}

\begin{proof}
    Let \(n=[\,G:H\,]\). Since \(H\trianglelefteq G\), the subgroup \(\rho^{-1}(\,H\,)\) is normal in \(\Gamma\). Hence \(p_H\colon S_H\to S\) is a regular covering with deck group \(\Gamma/\rho^{-1}(\,H\,)\cong G/H\). By Equation~\eqref{eq:volume-naturality}, \(\textnormal{vol}(\rho_H)=n\,\textnormal{vol}(\rho)\). Since \(\rho_H\) is geometrisable, its volume is positive, and therefore \(\textnormal{vol}(\rho)>0\). Moreover, \(d(\rho_H)=d(\rho)\). If \(d(\rho)>1\), Theorem~\ref{thm:discrete-geometrisation} implies that \(\rho\) is geometrisable. If \(d(\rho)=1\), then the geometrisability of the torsion-free representation \(\rho_H\) implies that \(\rho_H\) is faithful; since \(\ker\rho\subset\pi_1(S_H)\), it follows that the representation \(\rho\) is faithful as well. Thus Theorem~\ref{thm:discrete-geometrisation} again shows that \(\rho\) is geometrisable.

    Let \(\sigma\) be a branched hyperbolic structure on \(S\) with holonomy \(\rho\), and set \(\sigma_H=p_H^*\sigma\). Then \(\sigma_H\) has holonomy \(\rho_H\). Finally, for every deck transformation \(a\in G/H\), the identity \(p_H\circ a=p_H\) gives
    \begin{equation}
        a^*\sigma_H=a^*p_H^*\sigma=(p_H\circ a)^*\sigma=p_H^*\sigma=\sigma_H.
    \end{equation}
    Hence \(\sigma_H\) is invariant under the prescribed deck action of \(G/H\).
\end{proof}

In what follows, we shall define \textit{symmetric} those branched hyperbolic structures supplied by Theorem \ref{cor:invariant-cover}. Before moving on, we also notice that such a normal subgroup can always be obtained by replacing a given torsion-free subgroup \(K<G\) of finite index by its normal core
\begin{equation}
    H=\bigcap_{g\in G}gKg^{-1}.
\end{equation}
It is the kernel of the action of \(G\) on the finite coset set \(G/K\), so it is normal and has finite index. Moreover, it is torsion-free because it is contained in \(K\). Theorem \ref{cor:invariant-cover} asserts the existence of an invariant realisation. It does not require that a previously chosen realisation of \(\rho_H\) be preserved.

Theorem \ref{cor:invariant-cover} produces a branched hyperbolic structure on \(S_H\) which is invariant under the prescribed action of \(G/H\). We now describe how this action sits inside the full group of orientation-preserving isometries of the resulting structure. For this purpose, let \(\sigma_H\) be a branched hyperbolic structure on \(S_H\) with holonomy \(\rho_H\) and image \(H\). Since \(H\) is torsion-free and cocompact, the developing map descends to a branched covering
\begin{equation}
    f_H\colon S_H\longrightarrow \Sigma_H\defeq\HH^2/H.
\end{equation}
The metric \(\sigma_H\) is the pull-back by \(f_H\) of the hyperbolic metric on \(\Sigma_H\). Let \(\varphi\in\operatorname{Isom}^+(\sigma_H)\). Upon lifting \(\varphi\) to the universal cover, there exists an element \(g\in\pslr2\) such that
\begin{equation}
    \textnormal{dev}\circ\widetilde\varphi
    =
    g\circ\textnormal{dev}.
\end{equation}
Since the image of the holonomy is \(H\), the element \(g\) normalises \(H\). Notice that changing the lift of \(\varphi\) changes \(g\) by an element of \(H\). We therefore obtain a natural homomorphism
\begin{equation}\label{eq:isometry-normaliser-map}
    \Theta\colon
    \operatorname{Isom}^+(\sigma_H)
    \longrightarrow
    \frac{N_{\pslr2}(H)}{H}.
\end{equation}

The kernel of \(\Theta\) consists of those isometries \(\varphi\in\operatorname{Isom}^+(\sigma_H)\) satisfying the identity \(f_H\circ\varphi=f_H\). Thus
\begin{equation}
    \ker\Theta=\operatorname{Deck}(f_H)
\end{equation} the group of automorphisms of the branched covering \(f_H\colon S_H\to\Sigma_H\). Let \(N(f_H)\subset N_{\pslr2}(H)\) be the subgroup consisting of those elements whose induced isometry of \(\Sigma_H\) lifts through \(f_H\) to an isometry of \(\sigma_H\). We then obtain an exact sequence
\begin{equation}\label{eq:isometry-exact-sequence}
    1
    \longrightarrow
    \operatorname{Deck}(f_H)
    \longrightarrow
    \operatorname{Isom}^+(\sigma_H)
    \xrightarrow{\ \Theta\ }
    N(f_H)/H
    \longrightarrow
    1.
\end{equation}

For the realisation supplied by Theorem~\ref{cor:invariant-cover}, the prescribed deck action gives an inclusion
\begin{equation}
    G/H
    \hookrightarrow
    \operatorname{Isom}^+(\sigma_H).
\end{equation}
Under \(\Theta\), this subgroup is identified with the natural subgroup \(G/H\subset N_{\pslr2}(H)/H\). The exact sequence above therefore gives the following characterisation.

\begin{cor}\label{cor:full-symmetry-characterisation}
    Let \(\sigma_H\) be the branched hyperbolic structure supplied by Theorem~\ref{cor:invariant-cover}. Then the prescribed action of \(G/H\) is the full group of orientation-preserving isometries of \(\sigma_H\) if and only if
    \begin{equation}\label{eq:full-symmetry-condition}
        \operatorname{Deck}(f_H)=\{1\}
        \qquad\textnormal{and}\qquad
        N(f_H)=G.
    \end{equation}
    Equivalently, the natural inclusion \( G/H\hookrightarrow\operatorname{Isom}^+(\sigma_H)\) is an isomorphism provided the conditions in \eqref{eq:full-symmetry-condition} hold.
\end{cor}

\begin{proof}
    Suppose first that \(\operatorname{Isom}^+(\sigma_H)=G/H\). On the prescribed subgroup \(G/H\), the homomorphism \(\Theta\) is the natural inclusion into \(N_{\pslr2}(H)/H\), and is therefore injective. Hence
    \begin{equation}
        \operatorname{Deck}(f_H)=\ker\Theta=\{1\}.
    \end{equation}
    Moreover, the image of \(\Theta\) is \(G/H\). As a consequence of the short exact sequence \eqref{eq:isometry-exact-sequence}, the image of \(\Theta\) is also \(N(f_H)/H\), and therefore \(N(f_H)=G\). Conversely, suppose that \(\operatorname{Deck}(f_H)=\{1\}\) and that \(N(f_H)=G\). Then \eqref{eq:isometry-exact-sequence} reduces to \(\operatorname{Isom}^+(\sigma_H)\cong G/H\). Since \(G/H\) is already realised as the prescribed subgroup of \(\operatorname{Isom}^+(\sigma_H)\), it follows that
    \begin{equation}
        \operatorname{Isom}^+(\sigma_H)=G/H
    \end{equation}
    as desired.
\end{proof}


\smallskip

\section{Geometrisation in prescribed strata}\label{sec:strata}

We now refine the geometrisation problem by prescribing the branching orders of the branched hyperbolic structure. More precisely, we ask when a discrete representation can be realised by a branched hyperbolic structure belonging to a prescribed stratum. We first derive the numerical conditions imposed by the Euler number and by the degree of the associated branched covering. We then show that, when the image of the representation is torsion-free, these conditions are also sufficient.

\subsection{Strata}\label{subsec:strata}

Let \(\mathcal{BH}(S)\) denote the deformation space of marked branched hyperbolic structures on \(S\) compatible with its orientation. Recall that every \(\sigma\in\mathcal{BH}(S)\) determines a branching divisor, say \(D_\sigma\), defined as
\begin{equation}
    D_\sigma=\sum_{i=1}^n k_i p_i,
\end{equation}
where the \(p_i\)'s are distinct and \(k_i\geq1\) is the branching order at \(p_i\). Let \(b\) be a positive integer such that \(b<2g-2\), and let \(\mu=(k_1,\ldots,k_n)\), with \(1\leq k_1\leq\cdots\leq k_n\), be a partition of \(b\). We denote by \(\mathcal{BH}(\,\mu\,)\subset\mathcal{BH}(S)\) the stratum consisting of branched hyperbolic structures whose branching divisor has multiplicities prescribed by \(\mu\). The holonomy map, as defined in \eqref{eq:holonomy-map}, restricts to a map
\begin{equation}
    \hol(\mu)\colon
    \mathcal{BH}(\,\mu\,)
    \longrightarrow
    \operatorname{Hom}\bigl(\pi_1(S),\pslr2\bigr)/\pslr2.
\end{equation}
Our aim is to determine which discrete representations belong to the image of \(\hol(\,\mu\,)\). By Lemma~\ref{lem:euler-branching}, the image of each stratum is contained in a single Euler-number component. Indeed, every \(\sigma\in\mathcal{BH}(\,\mu\,)\) satisfies
\begin{equation}\label{eq:stratum-euler-number}
    e\big(\,\hol(\sigma)\,\big)=\chi(S)+b.
\end{equation}
Here \(S\) is oriented in such a way that branched hyperbolic structures compatible with the fixed orientation have holonomies with negative Euler number.

\smallskip

\subsection{Necessary conditions}\label{subsec:strata-necessary} In the present section we provide necessary conditions for a representation to arise as the holonomy of some branched hyperbolic structure with prescribed branched points. Let \(\rho\colon\pi_1(S)\to\pslr2\) have discrete image \(G\) and positive volume. By Lemma~\ref{lem:noncompact-euler-zero}, the quotient \(\mathcal O=G\backslash\HH^2\) is a compact hyperbolic orbifold. We set
\begin{equation}\label{eq:strata-total-branching}
    d(\rho)
    \defeq
    \frac{\textnormal{vol}(\rho)}
         {\textnormal{vol}(\mathcal O)}.
\end{equation}
As shown in Section~\ref{sec:discrete-geometrisation}, \(d(\rho)\) is a positive integer. The following gives the numerical obstructions to realising \(\rho\) in a prescribed stratum.

\begin{prop}\label{prop:strata-necessary}
    Let \(\mu=(k_1,\ldots,k_n)\) be a partition of an integer \(b\), with \(0<b<2g-2\). Suppose that
    \( \rho\) arises as the holonomy of some branched hyperbolic structure in \(\mathcal{BH}(\,\mu\,)\).
    Then the following conditions hold:
    \begin{itemize}
        \item[1.] \(\textnormal{vol}(\,\rho\,)=-2\pi e(\rho)=-2\pi(\,\chi(S)+b\,)>0\)
        \item[2.] Let \(\mathcal O=\HH^2/\rho(\pi_1(S))\). Then the inequality \begin{equation}\label{eq:strata-local-degree-condition}
            \textnormal{vol}(\rho)\geq\max_{1\leq i\leq n}\{k_i+1\}\,\textnormal{vol}(\mathcal O),
    \end{equation}
    holds.
    \end{itemize}
\end{prop}

\begin{proof}
    Let \(\sigma\in\mathcal{BH}(\mu)\) have holonomy \(\rho\), and let \(p_1,\ldots,p_n\) be its branch points, of orders \(k_1,\ldots,k_n\), respectively. The first assertion follows immediately from Lemma~\ref{lem:euler-branching}, since \(\kappa(\sigma)=b\). We now prove the second assertion. Let \(G=\rho(\pi_1(S))\), define \(\mathcal O=G\backslash\HH^2\) and let \(q\colon\HH^2\to|\mathcal O|\) denote the quotient map. Since the developing map \(\textnormal{dev}_\sigma\colon\widetilde S\longrightarrow\HH^2\) is \(\rho\)-equivariant, the composition \(q\circ\textnormal{dev}\) is invariant under the action of \(\pi_1(S)\), and hence descends to a map \(f_\sigma\colon S\longrightarrow|\mathcal O|\). Its degree is determined by the volume. Indeed, integrating the pull-back of the hyperbolic orbifold area form gives
    \begin{equation}\label{eq:strata-degree-volume}
        \textnormal{vol}(\rho)
        =
        \deg(f_\sigma)\,
        \textnormal{vol}(\mathcal O).
    \end{equation}

    Fix a branch point \(p_i\) of order \(k_i\), and let \(m_i\geq1\) be the orbifold order of \(f_\sigma(p_i)\), with \(m_i=1\) when \(f_\sigma(p_i)\) is a regular point of \(\mathcal O\). Choose a local complex coordinate \(z\) centred at \(p_i\) and an orbifold uniformising coordinate \(w\) centred at \(\textnormal{dev}(p_i)\). By the definition of the branching order, the developing map is locally equivalent to \(w=z^{k_i+1}\). If \(u\) is a coordinate on the underlying surface \(|\mathcal O|\), then locally \(u=w^{m_i}\). As a consequence, \(f_\sigma\) is locally equivalent to  \(u=z^{m_i(k_i+1)}\), and therefore the equality
    \begin{equation}
        \deg_{p_i}(f_\sigma)=m_i(k_i+1)
    \end{equation}

    holds. Since the local degree cannot exceed the global degree \(m_i(k_i+1)\leq\deg(f_\sigma)\). On the other hand, since \(m_i\geq1\), it follows that \( k_i+1\leq\deg(f_\sigma)\) for every \(i=1,\ldots,n\). Hence
    \begin{equation}\label{eq:max ineq}
        \max_{1\leq i\leq n}\{k_i+1\}
        \leq
        \deg(f_\sigma)
    \end{equation}
    holds. By combining \eqref{eq:max ineq} with \eqref{eq:strata-degree-volume}, the desired inequality readily follows, namely 
    \begin{equation}
        \textnormal{vol}(\rho)
        \geq
        \max_{1\leq i\leq n}\{k_i+1\}\,
        \textnormal{vol}(\mathcal O),
    \end{equation}
    as required.
\end{proof}


\subsection{The torsion-free case}\label{subsec:strata-torsion-free} In the present section we show that conditions of Proposition \ref{prop:strata-necessary} are also sufficient for purely hyperbolic representations. For this purpose, let \(\rho\) be a purely hyperbolic representation and let \(G=\rho(\pi_1(S))\) be torsion-free. We extend the notation of Lemma \ref{lem:F-surjective-degree} as follows: let \(\Sigma=G\backslash\HH^2\) and set
\begin{equation}\label{eq:degree-rho-surface}
    d(\,\rho\,)\defeq\frac{\textnormal{vol}(\,\rho\,)}{\textnormal{vol}(\Sigma)}
\end{equation}
which is the degree of the homotopy class determined by \(\rho\).

We shall use two classical results about branched coverings of surfaces. First, for a closed target surface different from the sphere and the projective plane, admissible branch data satisfying the Hurwitz congruence can be realised by a primitive branched covering; see \cite{BGKZ03, BGKZ04}. Here \emph{primitive} means that the induced homomorphism on fundamental groups is surjective. Second, Gabai and Kazez proved that primitive generic branched coverings of fixed degree are classified up to strong equivalence; see \cite{GabaiKazez1987}. Recall that a branched covering is \emph{generic} if all its ramification points are simple with pairwise distinct ramification values, and that strong equivalence allows a homeomorphism of the source together with a homeomorphism of the target isotopic to the identity.

\begin{thm}\label{thm:strata-torsion-free}
    Let \(\rho\colon\pi_1(S)\to\pslr2\) be a purely hyperbolic representation with positive volume, and let \(\mu=(k_1,\ldots,k_n)\) be a non-empty partition of \(0<b<2g-2\). Then \(\rho\) arises as the holonomy of some branched hyperbolic structure in \(\mathcal{BH}(\,\mu\,)\) if and only if 
    \begin{equation}\label{eq:strata-characterisation-total}
         e(\rho)=\chi(S)+b
    \end{equation}
    and
    \begin{equation}\label{eq:strata-characterisation-local}
        k_i+1\leq d(\,\rho\,)
    \end{equation}
    for every \(i=1,\ldots,n\).
\end{thm}

\begin{proof}
    Necessity follows from Proposition~\ref{prop:strata-necessary}, therefore it remains to prove the sufficiency. We argue as follows. For each \(i=1,\ldots,n\), consider the partition of \(d(\,\rho\,)\) given by
    \begin{equation}\label{eq:prescribed-partition-degree-d}
        A_i=\bigl(k_i+1,
        \underbrace{1,\ldots\ldots\ldots,1}_{d(\,\rho\,)-k_i-1\ \textnormal{terms}}\bigr).
    \end{equation}
    Notice that this is well defined by \eqref{eq:strata-characterisation-local}. Let \(\mathcal D_\mu=(A_1,\ldots,A_n)\) be the corresponding branch datum. According to \cite{EdmondsKulkarniStong1984}, its total defect equals \(\nu(\mathcal D_\mu)=b\). Set \(G=\rho(\pi_1(S))\) and \(\Sigma=G\backslash\HH^2\). Since \(G\) is torsion-free and cocompact, \(\Sigma\) is a closed hyperbolic surface. Moreover, by combining the definition of \(d(\,\rho\,)\) with Lemma \ref{lem:volume-euler}, the equality \(e(\,\rho\,)=d(\,\rho\,)\,\chi(\Sigma)\) holds. It follows from \eqref{eq:strata-characterisation-total} that
    \begin{equation}\label{eq:strata-hurwitz-defect}
        \nu(\mathcal D_\mu)=d(\,\rho\,)\,\chi(\Sigma)-\chi(S).
    \end{equation}
    In particular, the total defect is even. We can therefore apply the realisation theorem for primitive branched coverings, see \cite{BGKZ03, BGKZ04}, to the datum \(\mathcal D_\mu\). This yields a primitive branched covering, say \(g\colon S'\longrightarrow\Sigma\), of degree \(d(\,\rho\,)\) with one ramification point of order \(k_i\) over the \(i\)-th branch value and no other ramification. By the Riemann--Hurwitz formula and \eqref{eq:strata-hurwitz-defect},
    \begin{equation}
        \chi(S')=d(\,\rho\,)\,\chi(\Sigma)-\nu(\mathcal D_\mu)=\chi(S).
    \end{equation}
    Thus \(S'\) is orientation-preservingly homeomorphic to \(S\). After precomposing \(g\) with such a homeomorphism, we may regard it as a primitive branched covering \(g\colon S\to\Sigma\) whose branch points have orders \(k_1,\ldots,k_n\). 
    
    We notice that, in principle, there is no reason for the induced mapping \(g_*\) to coincide with the prescribed representation \(\rho\). it thus remains to recover the prescribed marking. By Theorem~\ref{thm:discrete-geometrisation}, the representation \(\rho\) is the holonomy of a branched hyperbolic structure on \(S\). Since the image \(G\) is torsion-free, the developing map of this structure descends to a branched covering \(f\colon S\longrightarrow\Sigma\) of degree \(d(\,\rho\,)\) such that \(f_*=\rho\), up to conjugation in \(G\), by design. In particular, \(f_*\) is surjective, and hence \(f\) is primitive. The branched coverings \(f\) and \(g\) have the same degree, but their ramification profiles need not agree. We therefore replace both by generic branched coverings. More precisely, by a local perturbation near each ramification point, every ramification point of order greater than one can be split into simple ramification points with pairwise distinct ramification values. In this way we obtain primitive generic branched coverings \(f_{\mathrm{gen}},g_{\mathrm{gen}}\colon S\to\Sigma\) in the sense of Gabai-Kazez, see \cite[pag. 1]{GabaiKazez1987}, such that
    \begin{equation}
        f_{\mathrm{gen}}\simeq f,   \qquad      g_{\mathrm{gen}}\simeq g.
    \end{equation}
    In particular, \((f_{\mathrm{gen}})_*=f_*=\rho\) up to conjugation. We now invoke the strong-equivalence theorem of Gabai--Kazez, see \cite{GabaiKazez1987}. Since \(f_{\mathrm{gen}}\) and \(g_{\mathrm{gen}}\) are primitive generic branched coverings of the same degree, there exist orientation-preserving homeomorphisms \(h\colon S\to S\) and \(a\colon\Sigma\to\Sigma\), with \(a\) isotopic to the identity, such that
    \begin{equation}\label{eq:strong-equivalence-strata}
        a\circ f_{\mathrm{gen}}\simeq g_{\mathrm{gen}}\circ h.
    \end{equation}
    Upon passing to fundamental groups and by using the homotopies above, we obtain
    \begin{equation}
        g_*\circ h_*=\rho
    \end{equation}
    up to conjugation. Indeed, \(a_*\) is the identity up to an inner automorphism because \(a\) is isotopic to the identity. Finally, we may notice that precomposition with \(h\) does not alter the ramification orders of \(g\). Hence the branched covering \(g\circ h\colon S\to\Sigma\) has branch points of orders \(k_1,\ldots,k_n\) and induces the representation \(\rho\), up to conjugation. Pulling back the hyperbolic metric \(h_\Sigma\) gives a branched hyperbolic structure \(\sigma=(g\circ h)^*h_\Sigma\) belonging to \(\mathcal{BH}(\,\mu\,)\) and having holonomy \([\rho]\). This proves the sufficiency.
\end{proof}

\smallskip

\subsection{The case with torsion}\label{ssec:torsion_case} We finally consider the case of discrete representations with torsion. The distinction between the torsion-free case considered in \S\ref{subsec:strata-torsion-free} can be understood in terms of markings. Suppose first that \(G<\pslr2\) is torsion-free and cocompact, and let \(\Sigma=G\backslash\HH^2\). If two epimorphisms, say \(\rho_1,\rho_2\colon\pi_1(S)\twoheadrightarrow G\) have the same negative Euler number, then they are induced by primitive maps \(S\longrightarrow\Sigma\) of the same degree. By the strong-equivalence theorem of Gabai--Kazez \cite{GabaiKazez1987} they differ, up to conjugation in \(G\), by precomposition with an orientation-preserving homeomorphism of \(S\). Equivalently, \([\rho_1]\) and \([\rho_2]\) lie in the same mapping class group orbit. Thus, in the torsion-free case, fixing the image \(G\) and the degree determines the representation up to the action of the mapping class group.

This scenario drastically changes when \(G\) contains torsion. In this case \(G\cong\pi_1^{\mathrm{orb}}(\mathcal O)\), and the argument above no longer provides control of the orbifold marking. As a consequence, realising a group \(G\) as the image of some holonomy representation in a prescribed stratum is weaker than realising a prescribed representation \(\rho\colon\pi_1(S)\twoheadrightarrow G\). As we shall see, the additional information can be expressed by passing to a normal torsion-free finite-index subgroup and requiring the corresponding realisation on the finite cover to be equivariant with respect to the prescribed deck action.

Let \(H\trianglelefteq G\) be a torsion-free normal subgroup of finite index \(N\), and set \(Q=G/H\). As in \S\ref{sec:sym structures}, let \(p_H\colon S_H\to S\) be the regular covering corresponding to \(\rho^{-1}(H)\), and write \(\rho_H=\rho|_{\pi_1(S_H)}\). The deck group of \(p_H\) is naturally identified with \(Q\). We also set \(\Sigma_H=H\backslash\HH^2\). Since \(H\trianglelefteq G\), then the quotient group \(Q\) acts naturally by orientation-preserving isometries on \(\Sigma_H\), and the quotient is the orbifold \(\mathcal O=G\backslash\HH^2\). Let \(\mu=(k_1,\ldots,k_n)\) be a non-empty partition. We denote by \(\mu_H\) the partition obtained by repeating every entry of \(\mu\) exactly \(N\) times; namely,
\begin{equation}\label{eq:lifted-partition}
    \mu_H=
    (\underbrace{k_1,\ldots,k_1}_{N\ \textnormal{times}},
    \ldots,
    \underbrace{k_n,\ldots,k_n}_{N\ \textnormal{times}})=(k_1^N,\ldots,k_n^N).
\end{equation}
Notice that this is precisely the branching profile obtained by pulling back a structure in \(\mathcal{BH}(\,\mu\,)\) through \(p_H\). Before stating the result, we observe that the numerical conditions are unchanged after passing to the finite cover. Indeed, by naturality of the Euler number and of the Euler characteristic,
\begin{equation}\label{eq:finite-cover-euler-strata}
    e(\,\rho_H\,)=N\,e(\,\rho\,),
    \qquad
    \chi(\,S_H\,)=N\,\chi(\,S\,),
\end{equation}
whereas
\begin{equation}\label{eq:finite-cover-degree-strata}
    d(\,\rho_H\,)
    =
    \frac{\textnormal{vol}(\,\rho_H\,)}
         {\textnormal{vol}(\Sigma_H)}
    =
    \frac{N\,\textnormal{vol}(\,\rho\,)}
         {N\,\textnormal{vol}(\mathcal O)}
    =
    d(\,\rho\,).
\end{equation}
Consequently, if
\begin{equation}\label{eq:torsion-strata-total}
    \sum_{i=1}^n k_i=e(\rho)-\chi(S) \quad \Longrightarrow \quad \sum_{k\in\mu_H}k=e(\rho_H)-\chi(S_H)
\end{equation}

Moreover, the inequalities \(k_i+1\leq d(\,\rho\,)\) are exactly the local degree inequalities for the lifted partition \(\mu_H\). Hence Theorem~\ref{thm:strata-torsion-free} guarantees the existence of a realisation of \(\rho_H\) in \(\mathcal{BH}(\mu_H)\). However, an additional issue arises, namely whether such a realisation can be chosen compatibly with the prescribed \(Q\)-action.

\begin{thm}\label{thm:strata-torsion}
    Let \(\rho\colon\pi_1(S)\to\pslr2\) have discrete cocompact image \(G\) with torsion and positive volume. Let \(H\trianglelefteq G\) be a torsion-free finite index subgroup, and let \(p_H\colon S_H\to S\) be the covering map corresponding to \(\rho^{-1}(\,H\,)\). Set
    \(\rho_H=\rho|_{\pi_1(S_H)}\), the group \(Q=G/H\) with \(N=[G:H]\) and \(\Sigma_H=H\backslash\HH^2\). Let \(b\) be a positive integer such that \(0<b<2g-2\), let \(\mu=(k_1,\ldots,k_n)\) be a non-empty partition of \(b\) and let \(\mu_H=(k_1^N,\ldots,k_n^N)\). Then \(\rho\) arises as the holonomy of some branched hyperbolic structure in \(\mathcal{BH}(\,\mu\,)\) if and only if the following conditions hold:
    \begin{itemize}
        \item[1.] \(e(\,\rho\,)=\chi(S)+b\),
        \item[2.] \(k_i+1\leq d(\rho)\) for every \(i=1,\ldots,n\),
        \item[3.] the representation \(\rho_H\) admits a realisation \(\sigma_H\in\mathcal{BH}(\,\mu_H\,)\) for which the prescribed deck action of \(Q\) on \(S_H\) is realised by orientation-preserving isometries and is compatible with the natural \(Q\)-action on \(\Sigma_H\). More precisely, if \(f_H\colon S_H\to\Sigma_H\) denotes the branched covering determined by \(\sigma_H\), then \(f_H\) can be chosen \(Q\)-equivariant, namely
        \begin{equation}\label{eq:Q-equivariant-strata}
            f_H(q\cdot x)=q\cdot f_H(x)
        \end{equation}
        for every \(q\in Q\) and \(x\in S_H\).
    \end{itemize}
\end{thm}

Although Theorem~\ref{cor:invariant-cover} guarantees the existence of a branched hyperbolic structure with prescribed holonomy which is invariant under a prescribed finite group action, it does not control the stratum of the resulting structure. In particular, it is not clear in general whether a given representation can be realised in a prescribed stratum while simultaneously preserving the prescribed symmetry. Condition~\textnormal{(3)} in Theorem~\ref{thm:strata-torsion} requires precisely this compatibility. As shown in the proof, this condition is also necessary for a realisation of \(\rho\) in the prescribed stratum.

\begin{rmk}\label{rmk:equivariant-splitting}
    The prescribed symmetry is preserved under splitting of branch points. More precisely, by using the notation of Theorem \ref{thm:strata-torsion}, let \(\sigma_H\) be a \(Q\)-invariant branched hyperbolic structure on \(S_H\) with holonomy \(\rho_H\), and suppose that one \(Q\)-orbit consists of branch points of order \(k\). Given a partition
    \begin{equation}
        k=k_1+\cdots+k_s,
    \end{equation}
    one may perform the usual local splitting surgery simultaneously along the whole \(Q\)-orbit. This produces a \(Q\)-invariant branched hyperbolic structure with the same holonomy, in which the original orbit is replaced by \(s\) \(Q\)-orbits of branch points of orders \(k_1,\ldots,k_s\), respectively. Therefore, whenever a \(Q\)-invariant realisation is available in a given stratum, any refinement of its branching partition can be achieved without losing the prescribed \(Q\)-symmetry. The converse operation is more delicate: it is not immediate that several \(Q\)-orbits of branch points can be collapsed equivariantly so as to produce an arbitrary prescribed partition. This is precisely the difficulty in removing Condition~\textnormal{(3)} from Theorem~\ref{thm:strata-torsion}.
\end{rmk}

\begin{rmk}\label{rmk:torsion-strata-obstruction}
    Conditions~\textnormal{(1)} and~\textnormal{(2)} already imply, by Theorem~\ref{thm:strata-torsion-free}, that \(\rho_H\) admits some realisation in \(\mathcal{BH}(\,\mu_H\,)\). Thus condition~\textnormal{(3)} does not concern the existence of a realisation on the torsion-free cover. It asks whether such a realisation can be chosen with the prescribed \(Q\)-symmetry. In this sense, the additional obstruction in the presence of torsion is an equivariance obstruction rather than a new numerical branching condition.
\end{rmk}

\begin{proof}
    Suppose first that \(\rho\) is realised by a branched hyperbolic structure \(\sigma\in\mathcal{BH}(\,\mu\,)\). Conditions~\textnormal{(1)} and~\textnormal{(2)} follow from Proposition~\ref{prop:strata-necessary}. Consider the pull-back structure \(\sigma_H=p_H^*\sigma\). Since the covering map \(p_H\) is an unbranched covering of degree \(N\), then every branch point of \(\sigma\) has exactly \(N\) preimages, with the same branching order. Hence \(\sigma_H\in\mathcal{BH}(\mu_H)\). Its holonomy is \(\rho_H\) by design. Moreover, for every deck transformation \(q\in Q\), we have \(p_H\circ q=p_H\), and therefore \(q\) preserves the branched hyperbolic structure \(\sigma_H\), that is,
    \begin{equation}
        q^*\sigma_H=\sigma_H.
    \end{equation}
    Thus the prescribed \(Q\)-action is realised by orientation-preserving isometries of \(\sigma_H\). Notice that the construction is compatible with the natural action of \(Q=G/H\) on \(\Sigma_H\). Equivalently, the associated branched covering \(f_H\colon S_H\to\Sigma_H\) satisfies \eqref{eq:Q-equivariant-strata}. This proves the necessity of~\textnormal{(3)}.

     Conversely, suppose that conditions~\textnormal{(1)}--\textnormal{(3)} hold. By~\textnormal{(3)}, the prescribed deck action of \(Q\) preserves the branched hyperbolic structure \(\sigma_H\). Since the action of \(Q\) on \(S_H\) is free, the branched hyperbolic metric \(\sigma_H\) descends to a branched hyperbolic structure \(\sigma\) on \(S_H/Q=S\). Every \(Q\)-orbit in \(S_H\) has cardinality \(N\), and isometries preserve branching orders. Since the branching profile of \(\sigma_H\) is that of \(\mu_H\), the quotient structure has branching profile \(\mu\). Thus \(\sigma\in\mathcal{BH}(\,\mu\,)\).

     It remains to identify the holonomy of the descended structure. Let \(\beta\colon\pi_1(S)\longrightarrow\pslr2\) denote this holonomy. Since the pull-back of the descended structure to \(S_H\) is precisely \(\sigma_H\), after an overall conjugation we may assume that
    \begin{equation}
        \beta|_{\pi_1(S_H)}=\rho_H.
    \end{equation}
    Let \(\gamma\in\pi_1(S)\) and \(\delta\in\pi_1(S_H)=\rho^{-1}(H)\). Since \(\pi_1(S_H)\) is normal in \(\pi_1(S)\), the element \(\gamma\delta\gamma^{-1}\) again belongs to \(\pi_1(S_H)\). Hence
    \begin{equation}
        \beta(\gamma)\rho_H(\delta)\beta(\gamma)^{-1}=\beta(\gamma\delta\gamma^{-1})=\rho(\gamma\delta\gamma^{-1})=\rho(\gamma)\rho_H(\delta)\rho(\gamma)^{-1}.
    \end{equation}
    Since \(\rho_H\) is surjective onto \(H\), it follows that \(\rho(\gamma)^{-1}\beta(\gamma)\) centralises \(H\). The group \(H\) is cocompact and non-elementary, and hence its centraliser in \(\pslr2\) is trivial. Therefore
    \begin{equation}
        \beta(\gamma)=\rho(\gamma)
    \end{equation}
    for every \(\gamma\in\pi_1(S)\). Thus the descended structure has holonomy \([\rho]\), and hence \([\rho]\in\hol(\mu)\bigl(\mathcal{BH}(\mu)\bigr)\) as desired.
\end{proof}

Theorem~\ref{thm:strata-torsion} therefore reduces the prescribed-stratum problem for representations with torsion to an equivariant realisation problem on a finite torsion-free cover. Notice that the symmetry required in condition~\textnormal{(3)} is weaker than requiring \(Q\) to be the full orientation-preserving isometry group of the structure \(\sigma_H\). 

\bibliographystyle{amsalpha}
\bibliography{discrgeo}

\end{document}